\input ./style/arxiv-general.cfg
\documentclass[MSNbibl,number,citesort,dvips]{arxbj}
\makeatletter
   \@ifpackageloaded{graphicx}{}{\usepackage{graphicx}}
\makeatother
\usepackage{euscript,multirow}

\volume{22}
\issue{4}
\pubyear{2016}
\firstpage{2177}
\lastpage{2208}
\doi{10.3150/15-BEJ725}
\docsubty{FLA}

\makeatletter
\newcommand{\rrvert}{\vert}
\newcommand{\rrVert}{\Vert}
\newcommand{\llvert}{\vert}
\newcommand{\llVert}{\Vert}
\renewcommand{\mid}{|}
\newcommand{\iint}{\int\!\!\!\int}
\newcommand{\iiint}{\int\!\!\!\int\!\!\!\int}
\newtheorem{theorem}{Theorem}
\newtheorem{lemma}[theorem]{Lemma}
\newremark{remark}{Remark}
\newcommand{\indep}{\ensuremath{\rlap{$\perp$}\mkern2mu{\perp}}}
\newcommand{\E}{\mathbb{E}}
\renewcommand{\P}{\mathbb{P}}
\newcommand{\R}{\mathbb{R}}
\newcommand{\N}{\mathbb{N}}
\newcommand{\w}[1]{{\widehat{#1}}}
\newcommand{\lr}{ \longrightarrow}
\newcommand{\var}{\operatorname{Var}}
\renewcommand{\r}{ \rightarrow}
\makeatother

\begin{document}
\begin{frontmatter}

\title{Integral approximation by kernel smoothing}
\runtitle{Integral approximation by kernel smoothing}

\begin{aug}
\author[A]{\inits{B.}\fnms{Bernard}~\snm{Delyon}\thanksref{A}\ead[label=e1]{bernard.delyon@univ-rennes1.fr}}
\and
\author[B]{\inits{F.}\fnms{Fran\c cois}~\snm{Portier}\corref{}\thanksref{B}\ead[label=e2]{francois.portier@gmail.com}}
\address[A]{Institut de recherches math\'ematiques de Rennes (IRMAR),
Campus de Beaulieu, Universit\'e de Rennes 1, 35042 Rennes C\'edex,
France. \printead{e1}}
\address[B]{Institut de Statistique, Biostatistique et Sciences
Actuarielles (ISBA), Universit\' e catholique de Louvain, Belgique.
\printead{e2}}
\end{aug}

%
\received{\smonth{9} \syear{2014}}
%
\revised{\smonth{3} \syear{2015}}

%
\begin{abstract}
Let $(X_1,\ldots,X_n)$ be an i.i.d. sequence of random variables in
$\mathbb R^d$, $d\geq 1$. We show that, for any function $\varphi
:\mathbb R^d\rightarrow\mathbb R$, under regularity conditions,
\[
n^{1/2} \Biggl(n^{-1} \sum_{i=1}^n
\frac{\varphi(X_i)}{\widehat f^(X_i)}-\int  \varphi(x)\,dx \Biggr)
\stackrel{
\mathbb P} {\longrightarrow} 0,
\]
where $\widehat f$ is the classical kernel estimator of
the density of $X_1$. This result is striking because it speeds up
traditional rates, in root $n$, derived from the central limit theorem
when $\widehat f=f$.
Although this paper highlights some applications, we mainly address
theoretical issues related to the later result.
We derive upper bounds for the rate of convergence in probability.
These bounds depend on the regularity of the functions $\varphi$ and
$f$, the dimension $d$ and the bandwidth of the kernel estimator
$\widehat f$.
Moreover, they are shown to be accurate since they are used as
renormalizing sequences in two central limit theorems each reflecting
different degrees of smoothness of $\varphi$.
As an application to regression modelling with random design, we
provide the asymptotic normality of the estimation of the linear
functionals of a regression function.
As a consequence of the above result, the asymptotic variance does not
depend on the regression function.
Finally, we debate the choice of the bandwidth for integral
approximation and we highlight the good behavior of our procedure
through simulations.
\end{abstract}

%
\end{frontmatter}

\section{Introduction}\label{s1}

Let $(X_1,\ldots,X_n)$ be an i.i.d. sequence of random variables in $\R
^d$, $d\geq 1$. We show that, for any function $\varphi:\R^d\r\R
$, under regularity conditions,
%
\begin{eqnarray}
\label{theresult} n^{-1} \sum_{i=1}^n
\frac{\varphi(X_i)}{\w f^{(i)}(X_i)}-\int  \varphi(x)\,dx = o_{\P}
\bigl(n^{-1/2}\bigr),
\end{eqnarray}
where $\w f^{(i)}$ is the classical leave-one-out kernel estimator of
the density of $X_1$ say $f$, defined by
\begin{eqnarray*}
\w f^{(i)}(x) = \bigl((n-1)h^d\bigr)^{-1}\sum
_{1\leq  j\leq  n, j\neq i} K\bigl(h^{-1} (x-X_j)
\bigr)\qquad\mbox{for every }x\in\R^d,
\end{eqnarray*}
where $K$ is a $d$-dimensional kernel and where $h$, called the
bandwidth, needs to be chosen and will certainly depend on $n$. Result
(\ref{theresult}) and the central limit theorem lead to the following
reasoning: when estimating the integral of a function that is evaluated
on a random grid $(X_i)$, whether $f$ is known or not, using a kernel
estimator of $f$ provides better convergence rates than using $f$ itself.

Result~(\ref{theresult}) certainly has some consequences in the field
of integral approximation. In this area, many deterministic as well as
random methods are available. Accuracy with respect to computational
time is the usual trade-off that allows to compare them. The advantages
of random over deterministic framework lie in their stability in
high-dimensional settings. For a comprehensive comparison between both
approaches, we refer to \cite{evans2000}.
Among random methods, \textit{importance sampling} is a widely used
technique that basically reduces the variance of the classical
Monte--Carlo integration through a good choice of the sampling
distribution $f$, called the sampler. Estimators are unbiased having
the form $n^{-1}\sum_{i=1}^n \varphi(X_i)/f(X_i)$ with $X_i\sim f$.
Regarding the mean squared error (MSE), the optimal sampler $f^*$ is
unique and depends on $\varphi$ (see Theorem 6.5 in \cite{evans2000},
page 176). Among others, parametric \cite{oh1992} and nonparametric
\cite{zhang1996} studies focused on the estimation of the optimal sampler.
Equation~(\ref{theresult}) indicates a new weighting of the
observations $\varphi(X_1),\ldots,\varphi(X_n)$.
Each weight $\w f^{(i)}(X_i)$ reflects how isolated is the point $X_i$
among the sample. Therefore, our estimator takes into account this
information by giving more weight to an isolated point. In summary our
procedure, which is adaptive to the design points
enjoys the following advantages:
\begin{itemize}
\item Faster than root $n$ rates,
\item one-step estimation based on a unique sample $(X_1,\ldots,X_n)$,
\item each $X_i$ drawn from $f$, possibly unknown.
\end{itemize}
To the best of our knowledge, when the design is not controlled, no
such rates have been obtained.

In many semiparametric problems, it has been an important issue to
construct root $n$ estimators, possibly efficient \cite{bickel1993},
that rely on a kernel estimator of the nuisance parameter. Among
others, it was addressed by Stone in \cite{stone1975} in the case of
the estimation of a location parameter, by Robinson in \cite
{robinson1988} in the \textit{partially linear regression model}, or by
H\"ardle and Stoker in \cite{hardle1989} studying the \textit{single
index model}.
The result in equation~(\ref{theresult}), which would be seen as a
superefficient estimator in the Le Cam's theory, cannot be linked
actually to this theory since the quantity of interest $\int  \varphi
(x)\,dx$ does not depend on the distribution of $X_1$. As a result, the
link between our work and the semiparametric literature relies mainly
on the plug-in strategy we employed, by substituting the density $f$ by
a kernel estimator.

In this paper, we propose a comprehensive study of the convergence
stated in equation~(\ref{theresult}). A~similar result was originally
stated by Vial in \cite{vial2003} (Chapter~7, equation (7.27)), as a
lemma in the context of the \textit{multiple index model}. To the best
of our knowledge, this type of asymptotic result has not been addressed
yet as a particular problem. Our theoretical aim is to extend result~(\ref{theresult}) by:
\begin{longlist}[(A)]
\item[(A)] Being more precise about the upper bounds: How does the dimension
$d$, the window $h$, the regularity of $\varphi$ and $f$, impact these bounds?
\item[(B)] Showing central limit theorems by specifying the regularity of
$\varphi$.
\end{longlist}
To achieve this program, we need to introduce a corrected version of
the estimate~(\ref{theresult}) for which the bias has been reduced.
First, the corrected estimator is shown to have better rates of
convergence than the initial one. Second,\vspace*{1pt} it is shown to be
asymptotically normal with rates $nh^{d/2}$ in the case where $\varphi$
is very regular, and with rates $ (nh^{-1})^{1/2}$ in a special case in
which $\varphi$ jumps at the boundary of its support. To compute the
asymptotic distribution, we rely on the paper by Hall \cite{hall1984},
where a central limit theorem for completely degenerate $U$-statistics
has been obtained. An important point is that we have succeeded in
proving our result with much weaker assumptions on the regularity of
$\varphi$ than on the regularity of $f$. For instance, equation~(\ref
{theresult}) may hold even when $\varphi$ has some jumps.
However, the estimation of $f$ is subject to the \textit{curse of
dimensionality}, that is, $f$ is required to be smooth enough regarding
the dimension of $X_1$.

Our aim is also to link equation~(\ref{theresult}) to nonparametric
regression with random design, that is, the model $Y_i=g(X_i)+\sigma
(X_i)e_i$ with $g$ unknown and $e_i$ i.i.d. with $e_i\,\indep\, X_i$.
In particular, we obtain the asymptotic normality for the estimators of
the linear functionals of $g$. Thanks to the fast rates detailed
previously, the asymptotic distribution does not depend on the function $g$.

The paper is organized as follows. Section~\ref{s2} deals with
technical issues related to equation~(\ref{theresult}). In particular,
we examine the rates of convergence of~(\ref{theresult}) according to
the choice of the bandwidth, the dimension and the regularity of the
functions $\varphi$ and $f$. Section~\ref{s3} is dedicated to the
convergence in distribution of our estimators. In Section~\ref{s4}, we
show how to apply equation~(\ref{theresult}) to the problem of the
estimation of the linear regression functionals. Finally, in
Section~\ref{s5}, we give some simulations that compare our method with
the traditional Monte--Carlo procedure for integration. The proofs and
the technicalities are postponed in Section~\ref{s6} at the end of the paper.

\section{Rates of convergences faster than root $n$}\label{s2}

In this section, we first provide upper bounds on the rates of
convergence in probability of our estimators. Our main purpose is to
show that rates faster than root $n$ hold in a wide range of parameter
settings for the estimation of $\int\varphi(x)\,dx$. Second, we argue
that those faster than root $n$ rates have no reason to hold when
estimating other functionals of the type $f\mapsto\int T(x,f(x)) \,dx$.

\subsection{Main result}\label{s21}

Let $Q\subset\R^d$ be the support of $\varphi$. The quantity $I(\varphi
)= \int \varphi(x)\,dx $ is estimated by
\begin{eqnarray*}
\w I(\varphi) = n^{-1} \sum_{i=1}^n
\frac{\varphi(X_i)}{\w f^{(i)}(X_i)}.
\end{eqnarray*}
Actually, this estimator can be modified in such a way that the leading
error term of its expansion vanishes asymptotically (see Remark~\ref
{cltnoncorrected} for more details).
For that, we define $\w v^{(i)}(x)$ as
\begin{eqnarray*}
\w v^{(i)}(x) = \bigl((n-1) (n-2)\bigr)^{-1} \sum
_{1\leq  j\leq  n, j\neq i} \bigl(h^{-d}K\bigl(h^{-1}(x-X_j)
\bigr)- \w f^{(i)}(x)\bigr)^2.
\end{eqnarray*}
It is, up to a factor $(n-1)^{-1}$, the leave-one-out estimator
of the variance of $h^{-d}K(h^{-1}(x-X_j))$.
The corrected estimator is
\begin{eqnarray*}
\w I_c(\varphi) = n^{-1} \sum
_{i=1}^n \frac{\varphi(X_i)}{\w f^{(i)}(X_i)} \biggl(1-
\frac{\w v^{(i)}(X_i)}{\w f^{(i)}(X_i)^2} \biggr).
\end{eqnarray*}
To state our main result about the convergences of $\w I(\varphi)$ and
$\w I_c(\varphi) $, we define the Nikolski class
of functions $\EuScript H(s,M)$ of regularity $s=k+\alpha$,
$k\in\mathbb N$, $0<\alpha\leq 1$, with constant $M>0$,
as the set of bounded and $k$ times differentiable functions
$\varphi$ whose all derivatives of order $k$ satisfy
\cite{tsybakov2009}
\begin{eqnarray*}
\int\bigl(\varphi^{(l)}(x+u)-\varphi^{(l)}(x)
\bigr)^2\,dx\leq  M\llvert u\rrvert^{2\alpha},\qquad l=(l_1,
\dots, l_d), \sum_{i=1}^d
l_i\leq  k,
\end{eqnarray*}
where $\llvert \cdot\rrvert $ stands for the Euclidean norm and the
$l_i$'s are natural integer.
Be careful that $k$ cannot be equal to $s$.
We say that $K$ is a kernel with order $r\in\N^*$ as soon as $K:\R
^d\mapsto\R$ is bounded and satisfies
\begin{eqnarray*}
\int K(x)\,dx &=& 1,\qquad\int x ^l K(x) \,dx =0,\qquad l=(l_1,
\dots, l_d), 0<\sum_{i=1}^d
l_i\leq  r-1 
\end{eqnarray*}
with the notation $x^l =x_1^{l_1}\times\cdots\times x_d^{l_d} $.
The following assumptions are needed to show our first result, they are
discussed after the statement.
\begin{enumerate}[(A1)]
\item[(A1)] For some $s>0$ and $M>0$, the support of $\varphi$
is a compact set $Q\subset\R^d$ and $\varphi$ is $\EuScript H(s,M)$ on
$\R^d$.
\item[(A2)] For some integer $r\geq 1 $, the variable $X_1$
has a bounded density $f$ on $\R^d$ such that its $r$th
order derivatives are bounded.
\item[(A3)] For every $x\in Q$, $f(x)\geq  b>0$.
\item[(A4)] The kernel $K$ has order $r$ and $\int K(x)\,dx =1$.
Moreover, there exists $ m_1>0$ and $m_2>0$ such that, for every $x\in
\R^d$, $\llvert K(x) \rrvert \leq  m_1\exp(-m_2\llvert x
\rrvert )$. In addition $K$ is symmetric: $K(x)=K(-x)$.
\end{enumerate}

The next theorem is proved in Section~\ref{s6}.


\begin{theorem}\label{thelemma}
Under the assumptions \textup{(A1)} to \textup{(A2)}, we have the following $O_\P$ estimates
\setcounter{equation}{0}
\renewcommand{\theequation}{\roman{equation}}
\begin{eqnarray}\label{borne1}
 n^{1/2} \bigl(\w I(\varphi) -I(\varphi) \bigr)
&=& O_\P\bigl( h^{s} + n^{1/2}h^{r} +  n^{-1/2}h^{-d} \bigr),
\\
n^{1/2} \bigl(\w I_c(\varphi)
 -I(\varphi)\bigr)&=& O_\P\bigl( h^{s} + n^{1/2}h^{r} +  n^{-1/2}h^{-d/2} +  n^{-1} h^{-3d/2} \bigr),\label{borne2}
\end{eqnarray}
which are valid if the sums inside the $O_\P$'s tend to zero.
\end{theorem}

\begin{remark}\label{rem1}
Assumption (A2) about the smoothness of $f$ is crucial to
guarantee a rate faster than root $n$ in Theorem~\ref{thelemma}. On the
one hand, one needs $r> d$ to obtain such a rate in equation~(\ref
{borne1}), on the other hand, $r> 3d/4$ suffices to get this rate in
equation~(\ref{borne2}).
Otherwise there does not exist $h$ such that the bounds in Theorem~\ref
{thelemma} go to $0$. This phenomenon is often referred as the \textit
{curse of dimensionality}.

In equation~(\ref{borne1}) (resp.,~(\ref{borne2})), when $h\propto
n^{-\gamma}$, the best choice of
$\gamma$ depends on $r$ and $s$; it balances two of the three (resp.,
four) terms while
letting the other one(s) smaller.
Precise rate acceleration for each situation is given in Table~\ref{tab1}.

\begin{table}
\tabcolsep=0pt
\caption{Best acceleration of convergence rate in Theorem \protect\ref{thelemma}.
Best rate acceleration $n^{-\beta}$ obtained with $h\propto n^{-\gamma} $}\label{tab1}
\begin{tabular*}{\tablewidth}{@{\extracolsep{\fill}}@{}lll@{}}
\hline
& $\beta$& $\gamma$ \\
\hline
Equation~(\ref{borne1}) \\
$2s\leq  r-d$  & $\frac{s}{2(s+d)} $    & $ {\frac{1}{2(s+d)}}$\\[5pt]
$0<r-d\leq 2s$ & $\frac{(r-d)}{2(r+d)} $& $ \frac{1}{r+d} $
\\[6pt]
Equation~(\ref{borne2}) &\\
$d\leq  r-d/2\leq 2s$ & $\frac{(r-d/2)}{2r+d} $ &$ \frac{1}{r+d/2} $ \\[3pt]
$d\leq 2s\leq  r-d/2$ & $ \frac{s}{2s+d} $& $ \frac{1}{2s+d} $ \\[5pt]
$r\leq 3d/2$  and $0<4r-3d\leq 6s$ & $\frac{4r-3d}{2(3d+2r)} $& $ \frac{3}{3d+2r} $ \\[5pt]
$2s\leq  d$ and $6s\leq 4r-3d$ & $ \frac{2s}{2s+3d} $&  $ \frac{2}{2s+3d} $ \\
\hline
\end{tabular*}
\end{table}

As in many semiparametric problems (see, e.g., \cite
{hardle1989}, Section~4.1), our estimator of $f$ is suboptimal with
respect to the density estimation problem (see \cite{stone1980}).
Indeed, to achieve the optimal rates in density estimation one would
need to take
$h \propto n^{-1/(2r+d)}$ which would even prevent $n^{1/2}h^r$ to go
to $0$ in Theorem~\ref{thelemma}.
A practical bandwidth selection is proposed Section~\ref{s5}.
\end{remark}

\begin{remark}\label{rem1bis}
Assumption (A2) prevents from bias problems in the estimation
of $f$ that may occur at the borders of $Q$. Indeed, if $f$ jumps at
the boundary of $Q$, then our estimate of $f$ would be asymptotically
biased and the rates provided in Theorem~\ref{thelemma} would not hold.
To get rid of this problem, if one knew the support of $f$, one could
correct by hand the estimator as, for instance, in \cite{jones1993}, or
might use Beta kernels as detailed in \cite{chen1999}.
\end{remark}

\begin{remark}\label{rem2}
Assumption (A3) basically says that $f$ is separated from $0$
on $Q$. The exponential bound on the kernel in assumption (A4)
guarantees that $f$ is estimated uniformly on $Q$ (see \cite
{devroye1980}). This helps to control the random denominators $\w
f^{(i)}(X_i)$'s in the expression of $\w I(\varphi)$ and $\w I_c(\varphi
)$. In the context of Monte--Carlo procedures for integral
approximation, assumptions (A2)~and~(A3)
are not that restrictive because one can draw the $X_i$'s from a
distribution smooth enough and whose support contains the integration domain.
\end{remark}

\begin{remark}\label{rem3}
The use of leave-one-out estimators $\w f^{(i)}$ and $\w v^{(i)}$ in
$\w I_c(\varphi)$ are not only justified by the simplification they
involve in the proofs. It also leads to better convergence rates.
Consider the term $R_0$ in the proof of equation~(\ref{borne2}) in
Theorem~\ref{thelemma}, when replacing the leave-one-out estimator of
$f$ by the classical one, $R_0$ remains a degenerate $U$-statistic but
with nonzero diagonal terms. It is possible to show that these terms
are leading terms of the resulting expansion. They imply a rate of
convergence of order $n^{-1/2}h^{-d}$ which is larger than the rate we
found for $\w I_c(\varphi)$.

However, concerning $\w I(\varphi)$, the leave-one-out estimator is not
necessary to get~(\ref{borne1}). The leave-one-out estimator being
indeed at a distance $O(n^{-1}h^{-d})$ from the ordinary one,
the change would made a difference of order at most $n^{-1/2}h^{-d}$ in
the left-hand side of~(\ref{borne1}),
which already appears in the right-hand side of~(\ref{borne1}).
\end{remark}

\begin{remark}\label{rem4}
The function class $\EuScript H(s,M)$ contains two interesting sets of
functions that provide different rates of convergence in Theorem~\ref
{thelemma}. First, if $\varphi$ is $\alpha$-H\"older on $\R^d$ with H\"
older constant $M_1$, and has bounded support, then $\varphi$
is $\EuScript H(\alpha,M_1)$ on $\R^d$. Second, if the support of
$\varphi$ is a convex body
(compact convex set with non-empty interior) and $\varphi$ is $\alpha
$-H\"older (with constant $M_1$)
\textit{inside its support}
(e.g., the indicator of a ball) then there exists $M_2>0$ such that
$\varphi$ is $\EuScript H(\min(\alpha,1/2),M_2)$ on $\R^d$ (see
Lemma~\ref{holderremark} in the Section~\ref{s6}). Then, because the
sum of two Nikolski functions is still Nikolski, the assumptions of
Theorem~\ref{thelemma} are valid for a wide range of integrand.
Moreover, note that a loss of
smoothness at the boundary of the support involves a loss in the rates
of convergence~(\ref{borne1}) and~(\ref{borne2}).
More precisely, whatever the smoothness degree of $\varphi$ inside its
support, if
continuity fails at the boundary, then the Nikolski regularity would be
at most $1/2$ and, therefore, the rates acceleration in Theorem~\ref
{thelemma} could not exceed $h^{1/2}$.
In Section~\ref{s3}, we study such an example and show a central limit
theorem with such a rate.
\end{remark}

\begin{remark}\label{rem4bis} The symmetry assumption in (A4) is
actually superfluous,
but simplifies the proof, because in this case we do not have to distinguish
the convolution with $K(x)$ and the convolution with $K(-x)$.
\end{remark}

\subsection{On the generalization of Theorem \texorpdfstring{\protect\ref{thelemma}}{1}}\label{generalizing}

In view of the intriguing convergence rates stated in Theorem
\ref{thelemma}, one may be curious to know the behavior of our
estimator when estimating more general functionals with the form
\begin{eqnarray*}
I_T= \int T\bigl(x,f(x)\bigr)\,dx,
\end{eqnarray*}
where $T:\R^{d}\times\R^+ \r\R$. Following the same approach as
previously, the estimator we consider is
\begin{eqnarray*}
\w I_{T} = n^{-1} \sum_{i=1}^n
\frac{T (X_i,\w f^{(i)}(X_i) )}{\w f^{(i)}(X_i)}.
\end{eqnarray*}
It turns out that $ T$ given by $(x,y)\mapsto\varphi(x)$ is the only
case for which the rates are faster than root~$n$. For other
functionals and a wide range of bandwidth, $\sqrt n (\w I_{T}- I_{T})$
converges to a normal distribution. In view of the negative aspect of
this result with respect to the statement of Theorem~\ref{thelemma}, we
provide an
informal calculation of the asymptotic law of $\sqrt n (\w
I_{T}-I_{T})$. We require that (A2) to (A4) hold and
that $n^{}h^{2r}\r0$ and $ n^{}h^{2d}\r+\infty$ (the latter
guarantees faster than root $n$ rates in equation~(\ref{borne1})).
If $y\mapsto T(x,y)$ has a bounded (uniformly in $x$) second-order derivative,
using a Taylor expansion with respect to the second coordinate of $T$
(the first-order derivative of $T$ with respect to the second
coordinate is further
denoted by $\partial_2 T$), we have
\begin{eqnarray*}
&&n^{1/2} (\w I_{T}-I_{T})
\\
&&\quad = n^{-1/2} \sum_{i=1}^n
\biggl( \frac{T(X_i,f(X_i))}{\w f^{(i)}(X_i)}- I_{T}+\frac{\partial_2
T(X_i,f(X_i))
( \w f^{(i)}(X_i)-f(X_i))}{\w f^{(i)}(X_i)} \biggr) +
\widetilde R_2,
\end{eqnarray*}
where $\widetilde R_2$ can be treated by standard techniques of kernel
estimation (see equations~(\ref{unifb}) and~(\ref{sansnom}) for
details), this gives that, with probability going to $1$,
\begin{eqnarray*}
\llvert\widetilde R_2\rrvert\leq  C n^{-1/2}\sum
_{i=1}^n \frac{(\w f^{(i)}(X_i)-f(X_i))^2}{\w f^{(i)}(X_i)}=
O_{\P} \bigl(n^{1/2}h^{2r}+ n^{-1/2}h^{-d}
\bigr)=o_{\P}(1),
\end{eqnarray*}
where $C>0$ does not depend on $n$ or $h$. Then we write
\begin{eqnarray*}
\sqrt n (\w I_{T}-I_{T}) =\widetilde R_0+
\widetilde R_1+\widetilde R_2,
\end{eqnarray*}
with
\begin{eqnarray*}
\widetilde R_0&=& n^{-1/2} \sum
_{i=1}^n \biggl(\frac{T(X_i,f(X_i))}{\w f^{(i)}(X_i)} -I_T
-\frac{\partial_2 T(X_i,f(X_i))
f(X_i)}{\w f^{(i)}(X_i)}+\int\partial_2 T\bigl(x,f(x)\bigr) f(x)\,dx
\biggr),
\\
\widetilde R_1&=& n^{-1/2} \sum
_{i=1}^n \biggl( \partial_2 T
\bigl(X_i,f(X_i)\bigr)-\int\partial_2 T
\bigl(x,f(x)\bigr) f(x)\,dx \biggr).
\end{eqnarray*}
If\vspace*{1pt} $x\mapsto T(x,f(x))$ and $x\mapsto\partial_2 T(x,f(x))f(x) $ are
Nikolski, applying Theorem~\ref{thelemma} gives that $\widetilde
R_0=o_{\P}(1)$. As a consequence $\sqrt n (\w I_{T}-I_{T})=o_{\P}(1)$
if and only if the variance of $\widetilde R_1$ is degenerate, that is
equivalent to
\begin{eqnarray*}
\partial_2 T\bigl(X_i,f(X_i)\bigr) = c
\qquad\mbox{a.s.}
\end{eqnarray*}
If we want this to be true for a reasonably large class of distribution
functions, it would imply
\begin{eqnarray*}
\partial_2 T(x,y) = c \qquad\mbox{for all } (x,y)\in
\R^d\times\R^+,
\end{eqnarray*}
for which the solutions have the form $T(x,y)= \varphi(x) + cy$.

\section{Central limit theorem}\label{s3}

In the previous section, we derived upper bounds on the convergence
rates in probability under fairly general conditions. In this section,
by being a little more specific about the regularity of~$\varphi$, we
are able to describe precisely the asymptotic distribution of $\w
I_c(\varphi)-I(\varphi)$.
Actually the approach is to decompose the latter quantity as a sum of a
$U$-statistic $U_n$
plus a martingale $M_n$ with respect to the filtration $\{X_1,\ldots,
X_n\}$, plus a bias term $B_n$ that is non-random (see the beginning of
Section~\ref{s62} for the definitions of $U_n$, $M_n$, $B_n$). Then
existing results about the asymptotic behavior of completely degenerate
$U$-statistics \cite{hall1984} and martingales \cite{hall1980} will help
to derive the asymptotic distribution. We shall consider two cases.
First, we present the case where $\varphi$ is smooth enough so that the
dominant term is $U_n$, and second we study an example where $\varphi$
is not continuous at the boundary of its support. As a consequence, the
dominant term is $M_n$.

For $\w I (\varphi)-I(\varphi)$, the situation is less interesting
since for most of the choice of $h$ a (non-random) bias term leads the
asymptotic decomposition (see Remark~\ref{cltnoncorrected}).

\subsection{Smooth case}

The smooth case corresponds to situations where the functions $f$ and
$\varphi$ are smooth enough, that is, $r>3d/2$ and $2s>d$. This is
highlighted by the assumptions on the bandwidth in the next theorem.

\begin{theorem}\label{tclphiregular}
Under the assumptions \textup{(A1)} to \textup{(A4)}, if $nh^{2d}\r
+\infty$, $nh^{r+d/2} \r0$ and $nh^{2s+d} \r0$, the random variable
$ nh^{d/2}(\w I_c(\varphi)-I(\varphi) )$ is asymptotically normally
distributed with zero-mean and variance given by\vspace*{-2pt}
\begin{eqnarray*}
\int\biggl(\int\bigl(K(u+v)-K(v)\bigr)K(u)\,du \biggr)^2\,dv \int
\varphi(x)^2f(x)^{-2} \,dx.
\end{eqnarray*}
\end{theorem}

The assumptions on the bandwidth are not satisfied by the optimal
bandwidths displayed in Table~\ref{tab1}. This is, in fact, a
presentation issue. Indeed we have chosen to make the bias term $B_n$
vanish so that any optimal bandwidth that balances the bias and the
variance is excluded. We could have proceeded the other way around, by
stating that $nh^{d/2}(\w I_c(\varphi)-I(\varphi) -B_n)$ has the same
limiting distribution as in Theorem~\ref{tclphiregular}, provided that
$ nh^{2d}\r+\infty$ and $nh^{2\min(r,s)+d} \r0$. One can verify that
this holds true for the optimal bandwidth given in the first line of
Table~\ref{tab1} for equation~(\ref{borne2}).

\subsection{A non-smooth example}

We are interested in the case where $\varphi$ is not sufficiently
regular so that $M_n$ is no longer negligible with respect to $U_n$,
that is, $nh^{2\min(r,s)+d}$ does not go to $0$.
This occurs whenever $s<d/2$.
In this case the variance is hard to compute since it depends on the
behavior of $M_n$
and therefore on the rate of convergence of the kernel regularization
of $\varphi$.
Hence, a precise description cannot be provided by considering usual regularity
classes, for example, H\"older, Nikolski or Sobolev since they only
provide bounds on the rate of
kernel regularization.
For this reason, we consider a particular case where the function
$\varphi$ is Nikolski inside
$Q$ and vanishes outside.
Typical functions we have in mind are the one that jump at the boundary
of their support.
Lemma~\ref{holderremark} informs us that such functions are Nikoslki
with regularity $1/2$. For $Q\subset\R^d$ compact and $x\in\partial
Q$, we define
\begin{eqnarray*}
L_Q(x)=\iint\min\bigl(\bigl\langle z,u(x)\bigr\rangle,\bigl\langle
z',u(x)\bigr\rangle\bigr)_+K(z)K\bigl(z'
\bigr)\,dz\,dz',
\end{eqnarray*}
where $u(x)$ is the unit normal outer vector of $Q$ at the point $x$.
We need the following assumption in place of (A1).
\begin{longlist}[(B1)]
\item[(B1)] For some $s>1/2$ and $M>0$, the support of $\varphi$
is a convex body $Q\subset\R^d$ with $\mathcal C^2$ boundary and
$\varphi$ is $\EuScript H(s,M)$ on $Q$.
\end{longlist}

\begin{theorem}\label{tclphinonregular}
Under\vspace*{1pt} the assumptions \textup{(A2)} to \textup{(A4)} and \textup{(B1)}, if $nh^{(3d+1)/2} \r +\infty$ and $nh^{2r-1} \r0$ the random
variable $ (nh^{-1})^{1/2} (\w I_c(\varphi)-I(\varphi)) $ is
asymptotically normally distributed with zero-mean and variance given by\vspace*{-2pt}
\begin{eqnarray*}
\int_{\partial Q} L_Q(x) \varphi(x)^2\,d
\mathcal H^{d-1}(x),
\end{eqnarray*}
where $\mathcal H^{d-1}$ stands for the $(d-1)$-dimensional Hausdorff measure.
\end{theorem}

\section{Application to nonparametric regression}\label{s4}

Equation~(\ref{theresult}) has applications in nonparametric
regression with random design. Let
%
\begin{eqnarray}
\label{modelad} Y_i= g(X_i) + \sigma(X_i)
e_i,
\end{eqnarray}
where $(e_i)$ is an i.i.d. sequence of real random variables with mean
$0$ and unit variance,
independent of the sequence $(X_i)$, and $\sigma:\R^d \r\R$ and $g:\R
^d \r\R$ are unknown functions. Let $Q\subset\R^d$ be a compact set
and $L_2(Q)$ be the Hilbert space of squared-integrable
functions on~$Q$.
Let $\psi\in L_2(Q)$ be extended to $\R^d$ by $0$ outside of $Q$ ($\psi
$ has compact support $Q$).
The inner product in $L_2(Q)$ between the regression function $g$ and
$\psi$, is given by
\begin{eqnarray*}
c = \int  g(x)\psi(x)\,dx,
\end{eqnarray*}
note that if $\psi$ belongs to a given basis of $L_2(Q)$, then $c$ is a
coordinate of $g$ in this basis. Among typical applications, we can
mention Fourier coefficients estimation for either nonparametric
estimation (see, e.g., \cite{hardle1990}, Section~3.3), or
location parameter estimation (see \cite{gamboa2007}). We also mention
the link with the estimation of the index in the \textit{single index
model} (see \cite{hardle1989}).

The estimation of the linear functionals of $g$ is a typical
semiparametric problem in the sense that it requires the nonparametric
estimation of the density $f$ of $X_1$
as a first step and then to use it in order to estimate a real parameter.
To the best of our knowledge, in the case of a regression with unknown
random design,
estimators that achieve root $n$ consistency have not been provided yet
(see, e.g., \cite{hardle1990} and the reference therein).
Our approach is based on kernel estimates $\w f^{(i)}$ of the density
of $X_1$ that are then
plugged into the classical empirical estimator of the quantity $\E[Y\psi
(X) f(X)^{-1}]$.
We define the estimator
\begin{eqnarray*}
\w c = n^{-1}\sum_{i=1}^n
\frac{Y_i \psi(X_i)}{\w f^{(i)}(X_i)},
\end{eqnarray*}
to derive the asymptotic of $\sqrt n (\w c -c )$, we use model~(\ref
{modelad}) to get the decomposition
\begin{eqnarray*}
\sqrt n (\w c -c )= A + B,
\end{eqnarray*}
with
\begin{eqnarray*}
A &=& n^{-1/2} \sum_{i=1}^n
\frac{\sigma(X_i)\psi(X_i)}{\w f^{(i)}(X_i)}e_i,
\\
B&=& n^{-1/2} \sum_{i=1}^n \biggl(
\frac{g(X_i)\psi(X_i)}{\w f^{(i)}(X_i)}-\int  g(x)\psi(x) \,dx \biggr).
\end{eqnarray*}
Roughly speaking, Theorem~\ref{thelemma} provides that $B$ is
negligible with respect to $A$. As a result, $A$ carries the weak
convergence of $\sqrt n (\w c -c )$ and, therefore, the limiting distribution
can
be obtained\vadjust{\goodbreak} making full use of the independence between the $X_i$'s and
the $e_i$'s.
In order to achieve such a program, this assumption is needed.
\begin{longlist}[(C1)]
\item[(C1)] For some $s>0$ and $M>0$, the support of $\psi$ is a
compact set
$Q\subset\R^d$ and both $\psi$ and $g$ are $\EuScript H(s,M)$ on $\R^d$.
\end{longlist}

The following theorem is proved in Section~\ref{s6}.

\begin{theorem}\label{thelemma2}
Under the assumptions \textup{(A2)} to \textup{(A4)}, and \textup{(C1)}, if $n^{1/2} h^r\r0$ and $n^{1/2}h^{d}\r
+\infty$, then the random variable $ n^{1/2}(\w c-c)$ is asymptotically
normally distributed with zero-mean and variance
\begin{eqnarray*}
v= \var\biggl(\frac{\sigma(X_1)\psi(X_1)}{f(X_1)} \biggr).
\end{eqnarray*}
\end{theorem}

\begin{remark}\label{rem5bis}
Let us compare $\w c$ with the appealing estimator
\begin{eqnarray*}
\widetilde c = n^{-1}\sum_{i=1}^n
\frac{Y_i \psi(X_i)}{f(X_i)}
\end{eqnarray*}
which requires the knowledge of $f$. First, if the signal is observed
without noise,
that is, $Y_i=g(X_i)$, then $n^{1/2}(\w c-c) $ goes to $0$ in
probability whereas
$\widetilde c$ is asymptotically normal. Secondly, when there is some
noise in the observed
signal, meaning that $\sigma(X_1) $ is not $0$,
the comparison can be made regarding their asymptotic variances. Since
we have
\begin{eqnarray*}
v\leq \var\bigl(n^{1/2}(\widetilde c-c)\bigr),
\end{eqnarray*}
it is asymptotically more efficient to plug the nonparametric estimator
of $f$ than to use $f$ directly.

\begin{remark}\label{rem5}
The set $Q$ reflects the domain where $g$ is studied. Obviously, the
more dense the $X_i$'s in $Q$, the more stable the estimation.
Nevertheless, it could happen that $f$ vanishes on some point on $Q$
and this is not taken into account by our framework. In such
situations, one may adapt the estimation from the sample by ignoring
the design points on which the estimated density takes too small
values. The estimator $\w c$ might be replaced by
\begin{eqnarray*}
n^{-1}\sum_{i=1}^n
\frac{Y_i \psi(X_i)}{\w f^{(i)}(X_i)} 1 _{\{\w f^{(i)}(X_i)>b\}},
\end{eqnarray*}
where $b>0$ will certainly depend on $n$. This method, often referred
as \textit{trimming}, has been employed in \cite{hardle1989} and \cite
{patilea2006} and guarantees computational stability as well as
theoretical properties. Even if such an approach is feasible here, it
seems far beyond the scope of the article.
\end{remark}
\end{remark}

\section{Simulations}\label{s5}

In this section, we provide some insights about the implementation and
the practical behavior of our integral approximation procedure. In
particular, we propose an adaptive procedure that selects the bandwidth
for the kernel smoothing.
While our theoretical study highlighted that our estimators suffers
from the \textit{curse of dimensionality} (see Remark~\ref{rem1}), our
simulation results confirm that the estimation accuracy of our methods
diminishes when the dimension increases. In dimension $1$, our
procedure outperforms by far the Monte--Carlo method. In moderate
sample size (from $200$ to $5000$) up to dimension $4$, our method
still realizes a significant improvement over the Monte--Carlo method.
The simulations are conducted under fairly general design distributions
that do not necessarily satisfy assumption (A2) (e.g., equation
(\ref{xunifo})).

\subsection{Kernel choice}
In the whole simulation study, our estimator of the density of the
design is based on the kernel
\begin{eqnarray*}
K(x) &=&\tfrac{1}2c_d^{-1}(d+1) \bigl(d+2-(d+3)
\llvert x \rrvert\bigr)1_{\llvert x \rrvert <1},
\\
c_d&=&\frac{2\pi^{d/2}}{d\Gamma(d/2)},
\end{eqnarray*}
where $c_d$ is the volume of the unit ball in dimension $d$. This
kernel is radial with order $3$.

\subsection{Bandwidth choice}
One may follow \cite{hardle1992} to select the optimal bandwidth
by a plug-in method.
It requires to optimize an asymptotic equivalent of the MSE with
respect to $h$. In Section~\ref{s3}, we highlighted that the limiting
distribution of $\w I(\varphi) - I(\varphi)$, and so the MSE, depends
heavily on the degree of smoothness of~$\varphi$. In practice, the
regularity of $\varphi$ is often unknown, as a result, we prefer a
simulation--validation
type strategy.

The idea is to pick the value $h$
which gives the best result for the estimation of the integral
$I(\widetilde\varphi)$ of a test function $\widetilde\varphi$ which
looks like $\varphi$,
and for which $I(\widetilde\varphi)$ is known.
We choose this test function as
%
\begin{eqnarray}
\label{sim1} \widetilde\varphi(x)=n^{-1}\sum
_{i=1}^n \frac{\varphi(X_i)}{\w f^{(i)}(X_i)}h_0^{-d}
\widetilde K \biggl(\frac{x-X_i}{h_0} \biggr),
\end{eqnarray}
where $\widetilde K$ is simply the Epanechnikov kernel
%
\begin{eqnarray}
\label{sim2} \widetilde K(x)&=&\tfrac{1}2c_d^{-1}(d+2)
\bigl(1-\llvert x \rrvert^2\bigr)1_{\llvert x \rrvert <1}.
\end{eqnarray}
Since we know that
\begin{eqnarray*}
I(\widetilde\varphi)=\int\widetilde\varphi(x)\,dx = n^{-1}\sum
_{i=1}^n \frac{\varphi(X_i)}{\w f^{(i)}(X_i)},
\end{eqnarray*}
we just take the value of $h$ for which the estimate $\w I(\widetilde
\varphi)$
is closest to $I(\widetilde\varphi)$; there is actually two values,
one for $\w I(\widetilde\varphi)$ and one for $\w I_c(\widetilde
\varphi)$.
The smoothing parameter $h_0$ is chosen using the rule of thumb given by
%
\begin{eqnarray}
\label{sim3} h_0=\sigma\biggl( \frac{d 2^{d+5}\Gamma(d/2+3)}{(2d+1)n}
\biggr)^{1/(4+d)},
\end{eqnarray}
where $\sigma^2$ is the mean of the estimated variances of each component
(see \cite{silverman1986}, Section~4.3.2).
The density estimates ${\w f^{(i)}(X_i)}$ in~(\ref{sim1}) are computed
with the same value $h_0$ and the same kernel.

We did not try to use a resampling method, thinking
that it is better to have $h$ adapted to the specific sample.

\subsection{First model}
\begin{figure}

\includegraphics{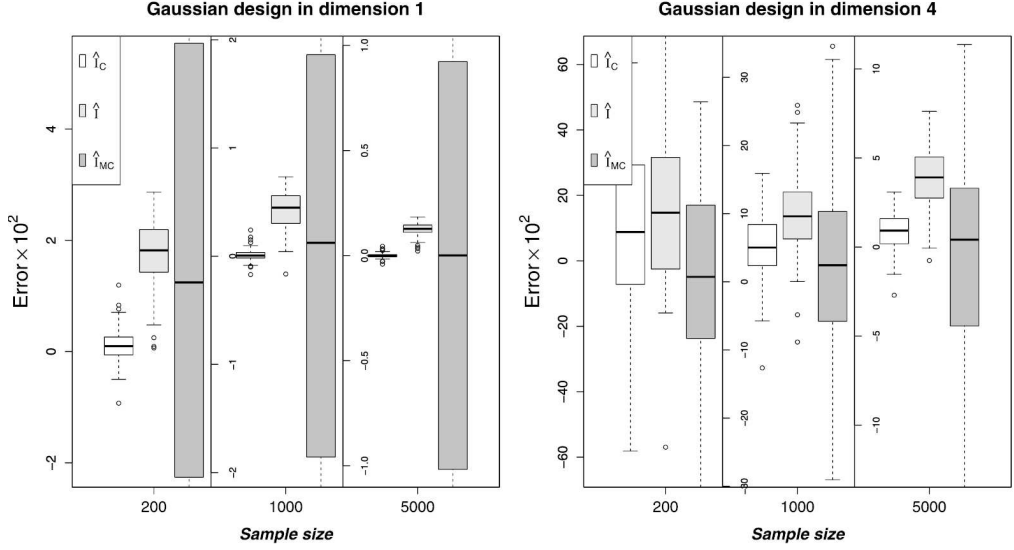}

\caption{Each boxplot is based on 100 estimates $\w I_c(\varphi)$, $\w
I(\varphi)$ and Monte--Carlo method noted $\w I_{\mathrm{MC}}$ for the first
model with different values of $n$ and $d$.}\label{simul1}
\end{figure}

In this model, $f$ is a normal distribution
\begin{eqnarray*}
X_i&\sim&\EuScript N\biggl(\frac{1}2,\frac{1}4 Id
\biggr),
\\
\varphi(x)&=&\prod_{k=1}^d2
\sin(\pi x_k)^21_{0\leq  x_k\le1}.
\end{eqnarray*}
The integral of $\varphi$ is 1. Figure~\ref{simul1} shows simulations
for different values of $n$ and $d$, and using equations~(\ref{sim1}),
(\ref{sim2}) and
(\ref{sim3}) for the choice of $h$.

%

\subsection{Second model}
In this second model, the assumptions are not satisfied
since the distribution is uniform over the unit cube, we have
%
\begin{eqnarray}
X_i&\sim&\EuScript U\bigl([0,1]^d\bigr),\label{xunifo}
\\
\varphi(x)&=&\prod_{k=1}^d2
\sin(\pi x_k)^21_{0\leq  x_k\le1}.
\end{eqnarray}
In spite of the fact that (A2) is not any more satisfied,
good results are still possible because $\varphi$ cancels at the boundary
of the cube.
For the choice of $h$, we used equation~(\ref{sim1}),
(\ref{sim2}) but, it is important to
constrain the function $\widetilde\varphi$ to have its support on the cube,
and a way to do this is to remove the boundary terms out of~(\ref
{sim1}) by choosing now
%
\begin{eqnarray}
\label{sim4} \widetilde\varphi(x)&=&\llvert J\rrvert^{-1}\sum
_{i\in J} \frac{\varphi(X_i)}{\w f^{(i)}(X_i)}h_0^{-d}
\widetilde K \biggl(\frac{x-X_i}{h_0} \biggr),
\nonumber\\[-8pt]\\[-8pt]\nonumber
J&=&\{i:h< X_{ij}<1-h,j=1\cdots d\}.
\nonumber
\end{eqnarray}

\begin{figure}

\includegraphics{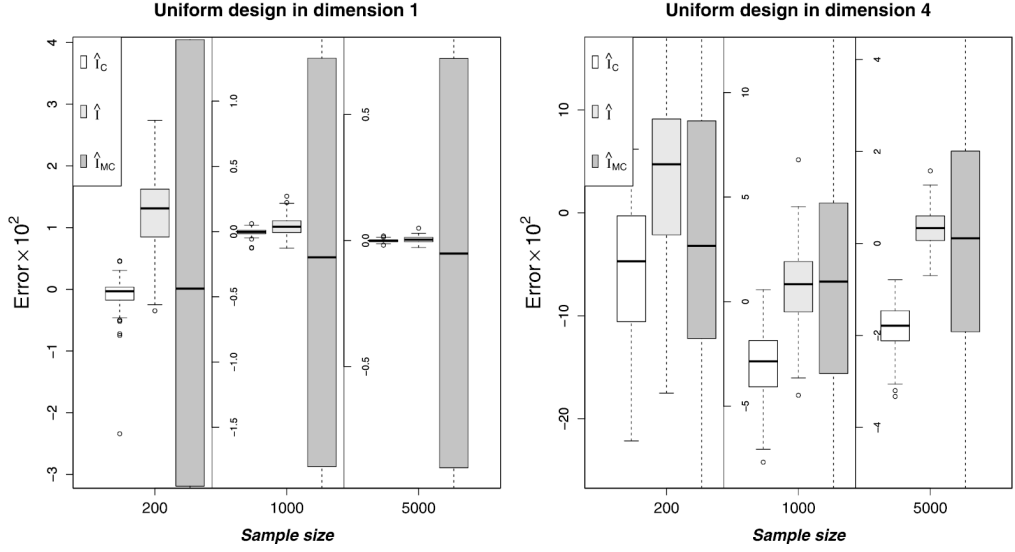}

\caption{Each boxplot is based on 100 estimates $\w I_c(\varphi)$, $\w
I(\varphi)$ and Monte--Carlo method noted $\w I_{\mathrm{MC}}$ for the second
model with different values of $n$ and $d$.}
\label{simul2}
\end{figure}

\noindent We could have done the other way around, use~(\ref{sim1}) and simulate uniformly
extra points at distance less than $h_0$ of the cube,
in order to cover the support of $\widetilde\varphi$. Figure~\ref
{simul2} shows the results of the simulations for different values of
$n$ and $d$ and using equations~(\ref{sim4}),~(\ref{sim2}) and
(\ref{sim3}) for the choice of~$h$.

\section{Proofs}\label{s6}

\subsection*{Notation}The Euclidean norm, the $L_p$ norm and the supremum
norm are, respectively, denoted by $\llvert \cdot\rrvert $, $\llVert
\cdot\rrVert _p$ and $\llVert \cdot\rrVert _\infty$. We introduce
$K_h(\cdot)= h^{-1}K(\cdot/h)$, and
\begin{eqnarray*}
K_{ij}&=&h^{-d}K\bigl(h^{-1}(X_i-X_j)
\bigr),
\\
\w f_i &=&\frac{1}{n-1} \sum_{1\leq  j\leq  n, j\neq i}
K_{ij},
\\
\w v_i &=& \frac{1}{(n-1)(n-2)} \sum_{1\leq  j\leq  n, j\neq i}
(K_{ij}- \w f_i)^2,
\end{eqnarray*}
and for any function $g:\R^d \r\R$, we define
\begin{eqnarray*}
g_{h}(x) = \int g(x-hu)K(u)\,du,
\end{eqnarray*}
and we put
\begin{eqnarray*}
\psi_q(x)&=&\frac{\varphi(x)}{f_h(x)^q},\qquad q\in\mathbb N,
\\
\widetilde\psi(x)&=& \biggl(\varphi(x)\frac{f(x)}{f_h(x)^2} \biggr)_h.
\end{eqnarray*}

\subsection{Proof of Theorem \texorpdfstring{\protect\ref{thelemma}}{1}}\label{s61}

We start by showing~(\ref{borne2}), then~(\ref{borne1}) will follow
straightforwardly.

\begin{pf*}{Proof of~(\ref{borne2})}
The following development reminiscent of the Taylor expansion
\begin{eqnarray*}
\frac{1}{\w f_i } = \frac{1} {f_h(X_i)} + \frac{f_h(X_i) -\w f_i }
{f_h(X_i)^2} +\frac{(f_h(X_i) -\w f_i )^2} {f_h(X_i)^3}
+\frac{(f_h(X_i)-\w f_i )^3} {\w f_i f_h(X_i)^3}, 
\end{eqnarray*}
allows us to expand our estimator as a sum of many terms, where the
density estimate $\w f_i$ is moved to the numerator,
with the exception of the fourth one. We will show that this last term
goes quickly to $0$. For the linearised terms, this is very messy
because the correct bound will be obtained by expanding also $\w f_i$
in those expressions.
In order to sort out these terms, we borrow from Vial \cite{vial2003}
the trick
of making appear a degenerate $U$-statistic in such a development (by
inserting the right quantity in $R_0$ below).
More explicitly, recalling that
\begin{eqnarray*}
\w I_c(\varphi) -I (\varphi) = n^{-1} \sum
_{i=1}^n \frac{\varphi(X_i)}{\w f_i} \biggl(1-
\frac{\w v_i}{\w f_i^2} \biggr) -I (\varphi),
\end{eqnarray*}
we obtain
%
\begin{eqnarray}
\label{decomplem1} \w I_c(\varphi) -I (\varphi)
=R_0+R_1+R_2+R_3+R_4+R_5,
\end{eqnarray}
with (we underbrace terms which have been deliberately introduced and removed)
\begin{eqnarray*}
R_0 &=& n^{-1} \sum_{i=1}^n
\bigl(\psi_1(X_i)-\psi_{2}(X_i)
{\w f_i } +\underbrace{\widetilde\psi(X_i)} -
\underbrace{\E\bigl[\psi_1(X_i)\bigr]} \bigr),
\\
R_1 &=& \int\bigl(\underbrace{f(x)f_h(x)^{-1}}-1
\bigr)\varphi(x) \,dx,
\\
R_2&=& n^{-1} \sum_{i=1}^n
\bigl(\psi_1(X_i)-\underbrace{\widetilde
\psi(X_i)} \bigr),
\\
R_3 &=& n^{-1} \sum_{i=1}^n
\psi_3(X_i) \bigl( \bigl( f_h(X_i)-
\w f_i \bigr)^2-\underbrace{\w v_i} \bigr),
\\
R_4 &=& n^{-1}\sum_{i=1}^n
\frac{\psi_3(X_i)
\w v_i}{\w f_i ^3 } \bigl(\underbrace{\w f_i ^3}-
f_h(X_i)^3 \bigr),
\\
R_5 &=& n^{-1} \sum_{i=1}^n
\psi_3(X_i) \frac{(f_h(X_i)-\w f_i )^3} {\w f_i },
\end{eqnarray*}
where $\w v_i$ appears to be a centering term in $R_3$.
We shall now compute bounds for each term separately.

%
\begin{longlist}[]
\item[\textit{Step} 1.]
$\llVert n^{1/2}R_0\rrVert _2=O(n^{-1/2} h^{-d/2})$.
Note that
\begin{eqnarray*}
R_0 &=& n^{-1} (n-1)^{-1}\sum
_{i\ne j} \bigl( \E\bigl[u_{ij}\llvert X_j
\bigr]-u_{ij}+E\bigl[u_{ij}\rrvert X_i
\bigr]-E[u_{ij}] \bigr),
\end{eqnarray*}
with $u_{ij}=\psi_{2}(X_i)K_{ij}$, is a degenerate $U$-statistic.
This is due to the fact that
\begin{eqnarray*}
\E[u_{ij}\mid X_i]&=&\psi_2(X_i)f_h(X_i)=
\psi_1(X_i),
\\
\E[u_{ij}\mid X_j]&=&(\psi_2f)_h(X_j)=
\widetilde\psi(X_j).
\end{eqnarray*}
The $n(n-1)$ terms in the sum are all orthogonal with $L_2$ norm
smaller than $\llVert u_{ij}\rrVert _2$, hence
\begin{eqnarray*}
(n-1)^2E\bigl[R_0^2\bigr]&\leq & \E
\bigl[u_{12}^2\bigr] \leq \llVert\psi_2
\rrVert_\infty^2E\bigl[K_{12}^2\bigr]
\leq  C_1 h^{-d},
\end{eqnarray*}
because of equation~(\ref{eq2}) in Lemma~\ref{moments}.\vspace*{1pt}

\item[\textit{Step} 2.]
$n^{1/2}R_1= O( n^{1/2}h^r)$. This is a consequence
of equation~(\ref{bochner2}) of Lemma~\ref{bochner},
and from assumption (A3).

\item[\textit{Step} 3.]
$\llVert n^{1/2}R_2\rrVert _2= O(
n^{1/2}h^r+h^{{s}})$. We can rearrange the function $\psi
_1(x)-\widetilde\psi(x)$
as
\begin{eqnarray*}
\psi_1(x)-\widetilde\psi(x)= \bigl(\psi_1(x)-
\psi_{1h}(x) \bigr)+ \bigl(\psi_{1h}(x)-\widetilde\psi(x)
\bigr),
\end{eqnarray*}
with
%
\begin{eqnarray}
\label{psitildeh} \bigl\llVert\psi_{1h}(x)-\widetilde\psi(x)\bigr
\rrVert
_\infty&=& \biggl\llVert\biggl(\psi_1(x)-\varphi(x)
\frac{f(x)}{f_h(x)^2} \biggr)_h\biggr\rrVert_\infty
\nonumber
\\
&\leq& \biggl\llVert\psi_1(x)-\varphi(x)\frac{f(x)}{f_h(x)^2}
\biggr\rrVert_\infty
\\
&=& \biggl\llVert\frac{\varphi}{f_h^2}(f_h-f)\biggr\rrVert
_\infty\leq  Ch^r,\nonumber
\end{eqnarray}
for some constant $C$, where the last inequality follows from equation
(\ref{bochner2}) in Lemma~\ref{bochner}. Then we have
\begin{eqnarray*}
R_2\leq \Biggl\llvert n^{-1}\sum
_{i=1}^n\bigl(\psi_{1h}(X_i)-
\psi_{1}(X_i)\bigr)\Biggr\rrvert+Ch^r,
\end{eqnarray*}
and by spliting the mean and the variance of the first term we get
\begin{eqnarray*}
\E\Biggl[ \Biggl(\frac{1}n \sum_{i=1}^n
\bigl(\psi_1(X_i) -\psi_{1h}(X_i)
\bigr) \Biggr)^2 \Biggr] &=&\E\bigl[\psi_1(X_1)
-\psi_{1h}(X_1)\bigr]^2+\frac{1}n\var
\bigl(\psi_1(X_1) -\psi_{1h}(X_1)
\bigr),
\end{eqnarray*}
and we conclude by equations~(\ref{norml1psih}) and~(\ref{norml2psih})
of Lemma~\ref{bochner} (it is an easy exercise to show that $\psi_1$ is
Nikolski with regularity $\min(r,s)$).

\item[\textit{Step} 4.]
$\llVert n^{1/2}R_3\rrVert _2=O(n^{-1/2}h^{-d/2})$.
We first express $R_3$ as a $U$-statistic. Set
\begin{eqnarray*}
U_i &=& \bigl(f_h(X_i)-\w f_i
\bigr)^2-\w v_i,
\end{eqnarray*}
and rewrite $R_3$ as
\begin{eqnarray*}
R_3 &=& n^{-1} \sum_{i=1}^n
\psi_3(X_i) U_i.
\end{eqnarray*}
Consider a sequence of real numbers $(x_j)_{1\leq  j\leq  p}$
and set
\begin{eqnarray*}
m&=&\frac{1}{p}\sum_{j=1}^px_j,
\\
v&=&\frac{1}{p(p-1)}\sum_{j=1}^p(x_j-m)^2=
\frac{1}{p(p-1)}\sum_{j=1}^p
\bigl(x_j^2-m^2\bigr),
\end{eqnarray*}
then
\begin{eqnarray*}
m^2-v= \biggl(1+\frac{1}{p-1} \biggr)m^2-
\frac{1}{p(p-1)}\sum_{j=1}^px_j^2=
\frac{2}{p(p-1)}\sum_{j< k}x_jx_k.
\end{eqnarray*}
Applying this with $x_j=K_{ij}-f_h(X_i)$ ($i$ is fixed) and $p=n-1$ we get
\begin{eqnarray*}
U_i&=&\frac{2}{(n-1)(n-2)}\sum_{j\ne i,k\ne i,j< k}
\bigl(K_{ij}-f_h(X_i) \bigr)
\bigl(K_{ik}-f_h(X_i) \bigr)
=\frac{2}{(n-1)(n-2)}\sum_{j< k}\xi_{ij}
\xi_{ik},
\end{eqnarray*}
with
\begin{eqnarray*}
\xi_{ij}&=&K_{ij}-f_h(X_i),
\\
\xi_{ii}&=&0.
\end{eqnarray*}
Then
\begin{eqnarray*}
R_3 &=& \frac{2}{n(n-1)(n-2)}\sum_i\sum
_{j< k}\psi_3(X_i)
\xi_{ij}\xi_{ik}.
\end{eqnarray*}
We are going to calculate $\E[R_3^2]$ by using the Efron--Stein
inequality (Theorem~\ref{efronstein})
and the moment inequalities~(\ref{eq1}) to~(\ref{eq3}) for
$\xi_{ij}$ stated in Lemma~\ref{moments}; in particular, by~(\ref
{eq1}), $\E[R_3^2]=\var(R_3)$.
Consider $R_3=f(X_1,\dots, X_n)$ as a function of the $X_i$'s and define
\begin{eqnarray*}
R'_3 &=& f\bigl(X'_1,X_2,
\ldots, X_n\bigr),
\\
\xi'_{1i}&=&h^{-d}K\bigl(h^{-1}
\bigl(X'_1-X_i\bigr)\bigr)-f_h
\bigl(X_1'\bigr),
\\
\xi'_{i1}&=&h^{-d}K\bigl(h^{-1}
\bigl(X_i-X'_1\bigr)\bigr)-f_h(X_i),
\\
\xi'_{ij}&=&\xi_{ij}\qquad\mbox{if }i\ne1
\mbox{ and }j\ne1,
\end{eqnarray*}
where $X'_1$ is a copy of $X_1$ independent from the sample $(X_1,\ldots
,X_n)$. Then by the Efron--Stein inequality (remember that $\xi_{ii}=0$)
\begin{eqnarray*}
\llVert R_3\rrVert_2&\leq & \biggl(
\frac{n }{ 2} \biggr)^{1/2} \bigl\llVert R_3-R'_3
\bigr\rrVert_2,
\end{eqnarray*}
which is of order
\begin{eqnarray*}
&&n^{-5/2}\biggl\llVert\sum_{j< k} \bigl(
\psi_3(X_1)\xi_{1j}\xi_{1k}-
\psi_3\bigl(X'_1\bigr)\xi_{1j}'
\xi_{1k}' \bigr) +\sum_i
\sum_{1< k}\psi_3(X_i) \bigl(
\xi_{i1}-\xi'_{i1}\bigr)\xi_{ik}
\biggr\rrVert_2
\\
&&\quad\leq  n^{-5/2} \biggl(\biggl\llVert\sum
_{j< k}\psi_3(X_1)\xi_{1j}
\xi_{1k}-\psi_3\bigl(X'_1\bigr)
\xi_{1j}'\xi_{1k}'\biggr\rrVert
_2 +\biggl\llVert\sum_{1< k}\sum
_i\psi_3(X_i) \bigl(
\xi_{i1}-\xi'_{i1}\bigr)\xi_{ik}
\biggr\rrVert_2 \biggr)
\\
&&\quad = n^{-5/2} \bigl(\llVert T_1\rrVert_2+
\llVert T_2\rrVert_2 \bigr).
\end{eqnarray*}
Noting that the terms in the first sum are orthogonal (by independence
of $\xi_{ij}$ and $\xi_{ik}$ conditionally to $X_i$ and~(\ref{eq1})) we obtain
\begin{eqnarray*}
\llVert T_1\rrVert_2&=& \frac{(n-1)^{1/2}(n-2)^{1/2}}{2^{1/2}} \bigl
\llVert
\psi_3(X_1)\xi_{12}\xi_{13}-
\psi_3\bigl(X'_1\bigr)\xi_{12}'
\xi_{13}'\bigr\rrVert_2
\\
&\leq& 2^{1/2} n \llVert\psi_3\rrVert_{\infty}
\llVert\xi_{12}\xi_{13}\rrVert_2
\\
&=&2^{1/2} n \llVert\psi_3\rrVert_{\infty}\E\bigl[
\E\bigl[\xi_{12}^2\xi_{13}^2\mid
X_1\bigr] \bigr]^{1/2}
\\
&=& 2^{1/2} n\llVert\psi_3\rrVert_{\infty}\bigl
\llVert\E\bigl[\xi_{12}^2\mid X_1\bigr] \bigr
\rrVert_2
\\
&=& O\bigl(nh^{-d}\bigr)
\end{eqnarray*}
by equation~(\ref{eq2}). Because the terms of the second sum
are orthogonal whenever the values of $k$ are different, we get
\begin{eqnarray*}
\llVert T_2\rrVert_2= (n-1)^{1/2}\biggl\llVert
\sum_i\psi_3(X_i) \bigl(
\xi_{i1}-\xi'_{i1}\bigr)\xi_{i2}
\biggr\rrVert_2.
\end{eqnarray*}
By first developing and then using that $X_1'$ is an independent copy
of $X_1$, we obtain
\begin{eqnarray*}
\biggl\llVert\sum_i\psi_3(X_i)
\bigl(\xi_{i1}-\xi'_{i1}\bigr)
\xi_{i2}\biggr\rrVert_2^2 &\leq&  n \E
\bigl[\psi_3(X_3)^2\bigl(\xi_{31}-
\xi'_{31}\bigr)^2\xi_{32}^2
\bigr]
\\
&&{}+n^2\bigl\llvert\E\bigl[\psi_3(X_3)
\psi_3(X_4) \bigl(\xi_{31}-
\xi'_{31}\bigr)\xi_{32}\bigl(
\xi_{41}-\xi'_{41}\bigr)\xi_{42}
\bigr]\bigr\rrvert
\\
& \leq& \llVert\psi_3\rrVert_{\infty}\bigl\{n \E\bigl[
\bigl(\xi_{31}-\xi'_{31}\bigr)^2
\xi_{32}^2 \bigr]
\\
&&{}+n^2 \E\bigl[\bigl\llvert\E\bigl[\bigl(\xi_{31}-
\xi'_{31}\bigr)\xi_{32}\bigl(
\xi_{41}-\xi'_{41}\bigr)\xi_{42}
\mid X_3,X_4\bigr]\bigr\rrvert\bigr]\bigr\}
\\
&=& \llVert\psi_3\rrVert_{\infty}\bigl\{ 2 n\E\bigl[
\xi_{31}^2\xi_{32}^2
\bigr]+2n^2 \E\bigl[\bigl\llvert\E[ \xi_{31}
\xi_{32}\xi_{41}\xi_{42}\mid X_3,X_4]
\bigr\rrvert\bigr]\bigr\}.
\end{eqnarray*}
Then by equation~(\ref{eq2}), we have $\E[\xi_{31}^2\xi_{32}^2 ]=\E[\E
[\xi_{31}^2\mid X_3]^2 ]\leq  C_1h^{-2d}$ and by equation~(\ref
{eq3}), we get
\begin{eqnarray*}
\E\bigl[\bigl\llvert\E[\xi_{31}\xi_{32}
\xi_{41}\xi_{42}\mid X_3,X_4]\bigr
\rrvert\bigr] &=&\E\bigl[\E[\xi_{31}\xi_{41}\mid
X_3,X_4]^2 \bigr]
\\
&\leq& 2\llVert f\rrVert_\infty^2 h^{-2d}\E
\bigl[\widetilde K\bigl(h^{-1}(X_4-X_3)
\bigr)^2 \bigr]+2\llVert f\rrVert_{\infty}^4
\\
&\leq & 2\llVert f\rrVert_\infty^3 h^{-d}\int
\widetilde K(u)^2\,d u+2\llVert f\rrVert_{\infty}^4,
\end{eqnarray*}
where $\widetilde K$ is defined in Lemma~\ref{moments}. Bringing
everything together and because $nh^d\r\infty$, it holds that
\begin{eqnarray*}
&& \bigl\llVert n^{1/2}R_3\bigr\rrVert_2
\leq  O\bigl( n^{-1}h^{-d}+ n^{-1/2}h^{-d/2}
\bigr)= O\bigl(n^{-1/2}h^{-d/2}\bigr).
\end{eqnarray*}

\item[\textit{Step} 5.]
$n^{1/2}R_4=O_{\P}(n^{-1} h^{-3d/2})$.
We start with a lower bound for $\w f_i$ by proving the existence of
$N(\omega)$ such that
%
\begin{eqnarray}
\forall n\geq  N(\omega), \forall i,\qquad\frac{b}{2}<\w
f_i<2\llVert f\rrVert_\infty.\label{unifb}
\end{eqnarray}
Notice that
\begin{eqnarray*}
\w f_i &=&\frac{n}{n-1} \biggl(\w f(X_i)-
\frac{1}{nh^d} K(0) \biggr),
\\
\w f(x)&=&\frac{1}{nh^d}\sum_{j=1}^nK
\bigl(h^{-d}(x-X_j)\bigr),
\end{eqnarray*}
due to the almost sure uniform convergence of $\w f$ to $f$ (Theorem~1
in \cite{devroye1980})
we have with probability~$1$ for $n$ large enough
\begin{eqnarray*}
\frac{2b}{3}<\inf_{x\in Q} \w f (x) \le\sup
_{x\in Q} \w f (x) < \frac{3}{2}\llVert f\rrVert
_\infty,
\end{eqnarray*}
and since assumption $nh^d\rightarrow\infty$,~(\ref{unifb}) follows.
We can now compute the expectation of $R_4$ restricted to $\{n\geq
N(\omega)\}$. Because $(a^3-b^3)=(a-b)(a^2+ab+b^2)$ for any real number
$a$ and $b$, and by the latter inequality, there exists a constant
$C>0$ which does not depend on $n$ or $h$, such that
\begin{eqnarray*}
\llvert R_4\rrvert1_{n>N(\omega)} &\leq&  C n^{-1}
\sum_{i=1}^n \bigl\llvert\w
f_i- f_h(X_i)\bigr\rrvert\w
v_i,
\end{eqnarray*}
we have by the Cauchy--Schwarz inequality
%
\begin{eqnarray}
\E\bigl[\llvert R_4\rrvert1_{n>N(\omega)}\bigr] &\leq&  C \E
\bigl[\bigl(\w f_1- f_h(X_1)
\bigr)^2\bigr]^{1/2}\E\bigl[\w v_1^2
\bigr]^{1/2}.\label{st50}
\end{eqnarray}
Applying the fact that for any real number $a$, $ \frac{1} p \sum
_{j=1}^p(x_j-\overline x)^2\leq \frac{1} p \sum_{i=1}^p(x_j-a)^2$
to $x_j=K_{1j}$, $p=n-1$ and $a=f_h(X_1)$, we obtain that
\begin{eqnarray*}
&&\w v_1 \leq \frac{1}{(n-1)(n-2)} \sum
_{j=2}^n \xi_{1j}^2,
\end{eqnarray*}
then using~(\ref{eq2})
%
\begin{eqnarray}\label{st51}
\nonumber
\E\bigl[\w v_1^2\bigr] &\leq&
(n-1)^{-1}(n-2)^{-2}\E\bigl[\xi_{12}^4
\bigr]+(n-1)^{-1}(n-2)^{-1}\E\bigl[\xi_{12}^2
\xi_{13}^2\bigr]
\\
&\leq&  C_1n^{-3}h^{-3d}+C_1n^{-2}h^{-2d}
\\
&\leq&  O\bigl(n^{-2}h^{-2d}\bigr),\nonumber
\end{eqnarray}
because $nh^d$ goes to infinity. On the other hand using equation~(\ref
{eq2}) again,
%
\begin{eqnarray}
\label{st52} \E\bigl[\bigl(\w f_1- f_h(X_1)
\bigr)^2 \bigr]&=&\frac{1}{n-1}\E\bigl[\xi_{1i}^2
\bigr] =O \bigl(n^{-1}h^{-d}\bigr).
\end{eqnarray}
Putting together~(\ref{st50}),~(\ref{st51}) and~(\ref{st52}),
\begin{eqnarray*}
\E\bigl[\llvert R_4\rrvert1_{n>N(\omega)}\bigr] &=& O\bigl(
n^{-1}h^{-d}n^{-1/2}h^{-d/2}\bigr)=O\bigl(
n^{-3/2}h^{-3d/2}\bigr).
\end{eqnarray*}
In particular by Markov's inequality
\begin{eqnarray*}
\P\bigl(n^{3/2} h^{3d/2}\llvert R_4\rrvert>A\bigr)
&\leq& \P\bigl(n^{3/2} h^{3d/2}\llvert R_4
\rrvert1_{n>N(\omega)}>A \bigr)+\P\bigl(n\leq  N(\omega) \bigr)
\\
&=& A^{-1}O(1)+\P\bigl(n\leq  N(\omega) \bigr).
\end{eqnarray*}
This proves the boundedness in probability of $n^{3/2} h^{3d/2}\llvert
R_4\rrvert $.\vspace*{2pt}

\item[\textit{Step} 6.]
$n^{1/2}R_5= O_\P( n^{-1}h^{-3d/2}+n^{-3/2} h^{-2d})$.
Following~(\ref{unifb}) since
\begin{eqnarray*}
\llvert R_5\rrvert1_{n>N(\omega)} \leq 2b^{-3}
\llVert\varphi\rrVert_\infty n^{-1} \sum
_{i=1}^n{\bigl\llvert\w f_i
-f_h(X_i)\bigr\rrvert^3},
\end{eqnarray*}
we can show the convergence in probability of the right-hand side term as in
Step 5. We have
indeed by the Rosenthal's inequality\footnote{\label{foot}For a
martingale $(S_i,\mathcal{F}_i)_{i\in\N}$
and $2\leq  p< +\infty$, we have $\E[\llvert S_n\rrvert
^p]\leq  C_2\{\E[(\sum_{i=1}^n\E[X_i^2\llvert \mathcal
{F}_{i-1}])^{p/2}]+\sum_{i=1}^n \E\rrvert X_i\mid^p\}$, where $X_i =
S_i -S_{i-1}$ (see, e.g., \cite{hall1980}, pp. 23--24).}
%
\begin{eqnarray}\label{rose}
\nonumber
\E\Biggl[n^{-1} \sum_{i=1}^n
\bigl\llvert\w f_i -f_h(X_i)\bigr\rrvert
^p \Biggr]&=&(n-1)^{-p}\E\Biggl[\Biggl\llvert\sum
_{i=2}^n \xi_{1i} \Biggr\rrvert
^p\Biggr]\nonumber
\\
&\leq&  C_2n^{-p}\bigl\{\bigl(n\E\bigl[
\xi_{12}^2\bigr]\bigr)^{p/2}+n\E\bigl[\llvert
\xi_{12}\rrvert^p\bigr]\bigr\}
\\
&\leq&  C_1C_2\bigl\{n^{-p/2}h^{-pd/2}+n^{1-p}h^{-(p-1)d}
\bigr\},\nonumber
\end{eqnarray}
where the latter inequality is due to equation~(\ref{eq2}). Hence, with $p=3$
\begin{eqnarray*}
\E\bigl[\llvert R_5\rrvert1_{n>N(\omega)} \bigr] \leq
C_1C_2\bigl\{n^{-3/2}h^{-3d/2}+n^{-2}h^{-2d}
\bigr\}
\end{eqnarray*}
and we conclude as in Step 5.

Putting together the steps 1 to 6, and taking into account, concerning $R_5$,
that $n^{-3/2} h^{-2d}=(n^{-1/2} h^{-d/2})(n^{-1} h^{-3d/2})$, we
obtain~(\ref{borne2}).
\end{longlist}

\begin{pf*}{Proof of~(\ref{borne1})}
For~(\ref{borne1}), we use a
shorter expansion which leads to an actually much simpler proof:
\begin{eqnarray*}
\frac{1}{\w f_i } = \frac{1} {f_h(X_i)} + \frac{f_h(X_i) -\w f_i }
{f_h(X_i)^2} +
\frac{(f_h(X_i) -\w f_i )^2} {\w f_i f_h(X_i)^2},
\end{eqnarray*}
and
\begin{eqnarray*}
\w I(\varphi) -I(\varphi) =R_0+R_1+R_2+R'_5,
\end{eqnarray*}
with
\begin{eqnarray*}
R'_5 &=& n^{-1} \sum
_{i=1}^n\psi_2(X_i)
\frac{(f_h(X_i)-\w f_i )^2} {\w f_i }.
\end{eqnarray*}
The terms $R_0$, $R_1$ and $R_2$ have already been treated in the steps
1, 2 and 3 of the proof of~(\ref{borne1}). The term $R'_5$ is bounded
exactly as $R_5$ but now we use~(\ref{rose}) with
$p=2$ instead of $p=3$, to obtain
\begin{eqnarray*}
\E\bigl[\bigl\llvert R'_5\bigr\rrvert1_{n>N(\omega)}\bigr)
\leq  C_1C_2n^{-1} h^{-d}
\end{eqnarray*}
and we get $n^{1/2}\llvert R_5'\rrvert =O_\P(n^{-1/2} h^{-d})$.
\end{pf*}

\subsection{Proofs of Theorems \texorpdfstring{\protect\ref{tclphiregular}}{2} and \texorpdfstring{\protect
\ref{tclphinonregular}}{3}}\label{s62}

Let us define
\begin{eqnarray*}
M_n &=& n^{-1}\sum_{i=1}^n
\psi_1(X_i)-{\widetilde\psi(X_i)} -\E
\bigl[\psi_1(X_1)-\widetilde\psi(X_1)
\bigr],
\\
U_n & =& n^{-1}(n-1)^{-1} \sum
_{i\neq j} c_{ij},
\\
B_{n}& =& \E\bigl[\psi_1(X_1)-\widetilde
\psi(X_1) \bigr]+\int\bigl(f(x)f_h(x)^{-1}-1
\bigr)\varphi(x) \,dx
\end{eqnarray*}
with $c_{jk}=a_{jk}-b_{jk}$, and for $j\neq k$,
\begin{eqnarray*}
a_{jk}& =& \E\bigl[\psi_3(X_1) {
\xi_{1j}\xi_{1k} } \mid X_j X_k
\bigr],
\\
b_{jk}&=&u_{jk} - \E\bigl[u_{jk}\llvert
X_j\bigr]-\E\bigl[u_{jk}\rrvert X_k\bigr]+
\E[u_{jk}],
\end{eqnarray*}
where $u_{jk}$ has been defined at the beginning of step 3. Both proofs
of Theorems~\ref{tclphiregular} and~\ref{tclphinonregular} rely on the
following lemma which turns Theorem~\ref{thelemma} in a suitable way
for weak convergence issues.

\begin{lemma}\label{decompustatmart}
Under the assumptions of Theorem \ref{thelemma}, we have
\begin{eqnarray*}
\w I_c(\varphi)-I(\varphi) = B_n+ U_n +
M_n +O_\P\bigl( n^{-3/2}h^{-3d/2}\bigr).
\end{eqnarray*}
Moreover, we have $B_n =O_\P(h^r) $, $U_n=O_{\P}(n^{-1}h^{-d/2}
)$ and $M_n=O_{\P}(n^{-1/2}(h^{s}+ h^{r}))$.
\end{lemma}

\begin{pf}
By using the decomposition~(\ref{decomplem1}) and since
$B_n+M_n=R_1+R_2$, we have
\begin{eqnarray*}
\w I_c(\varphi)-I (\varphi) &=&R_0+R_1+R_2+R_3+R_4+R_5
\\
&=& B_n+M_n+U_n+ (R_0+R_3-U_n)
+R_4+R_5.
\end{eqnarray*}
We have already shown that $R_4 +R_5= O_\P( n^{-3/2}h^{-3d/2}+n^{-2}
h^{-2d})$ (this is exactly steps 5~and~6 of the proof of Theorem~\ref{thelemma}). By
definition of $U_n$, we have
\begin{eqnarray*}
R_0+R_3-U_n = n^{-1}(n-1)^{-1}(n-2)^{-1}
\sum_{i}\sum_{j\neq k}
\bigl(\psi_3(X_i) {\xi_{ij}\xi_{ik}
}-a_{jk} \bigr)
\end{eqnarray*}
which is a completely degenerate $U$-statistic ($R_3$ is near to be
completely degenerate and $a_{jk}= \E[\psi_3(X_1) {\xi_{1j}\xi_{1k} }
\mid X_j X_k]$ appears as the good centering term). The order $2$
moments of this quantity are of order $ n^{-3}\E[\psi_3(X_1)^2 \xi
_{12}^2\xi_{13}^2 ]\propto n^{-3}h^{-2d} $. Hence, we have shown that
$R_0+R_3 - U_n=O_{\P}(n^{-3/2}h^{-d} )$, which completes the first part
of the proof. To obtain the bounds in probability, for $U_n$ we just
use step 1 and 4 of the proof of Theorem~\ref{thelemma}, for $M_n$ we
compute the $L_2$ norm as follows. We have
\begin{eqnarray*}
\llVert M_n\rrVert_2 &=& n^{-1/2} \bigl\llVert
\psi_1(X_1)-{\widetilde\psi(X_1)}\bigr
\rrVert_2
\\
&\leq&  n^{-1/2} \bigl(\bigl\llVert\psi_1(X_1)-
\psi_{1h}(X_1)\bigr\rrVert_2+\bigl\llVert
\psi_{1h}(X_1)-{\widetilde\psi(X_1)}\bigr
\rrVert_2\bigr)
\\
& \leq&  Cn^{-1/2}\bigl(h^{s}+h^{r}\bigr),
\end{eqnarray*}
for some constant $C$, where the last inequality is obtained using
equation~(\ref{psitildeh}) for the term in the right and equation~(\ref
{norml2psih}) in Lemma~\ref{bochner} for the term in the right.
\end{pf}

\begin{remark}\label{cltnoncorrected}
Under the assumption of Theorem~\ref{thelemma}, one may show that
\begin{eqnarray*}
\w I (\varphi)-I(\varphi)=\w I_{c}(\varphi)-I(
\varphi)+n^{-1} (n-1)^{-2} \sum_{i,j}
\psi_3(X_i) {\xi_{ij}^2
}+O_{\P}\bigl(\bigl(nh^{d}\bigr)^{-3/2}\bigr),
\end{eqnarray*}
where the $O_\P$ comes from $R_4$ and the other remainder term
corresponds to the diagonal term of the $U$-statistic $R_3$. This term
equals $(n-1)^{-1} \E[ \psi_{2} (X_1) (K_{12}-f_h(X_1))^2]=O
(n^{-1}h^{-d})$ plus $o_{\P} (n^{-1}h^{-d/2})$,
as a consequence, when $h$ is such that $nh^{2(s+d)}\r0$ and
$nh^{r+d}\r0$, the leading term of the decomposition is a constant.
\end{remark}

\subsubsection{Proof of Theorem \texorpdfstring{\protect\ref{tclphiregular}}{2}}
By Lemma~\ref{decompustatmart} and the assumptions on $h$ we have
\begin{eqnarray*}
n^{}h^{d/2}(B_n+M_n+R_4+R_5)
&=& O_\P\bigl( n^{3/2}h^{r+d/2}+n^{1/2}
\bigl(h^{r+d/2}+h^{s+d/2}\bigr) +n^{-1/2}h^{-3d/2}
\bigr)
\\
&=& o_{\P}(1).
\end{eqnarray*}
To derive the limiting distribution of $nh^{d/2}U_n$, we apply Theorem
1 in \cite{hall1984},
quoted below (Theorem~\ref{ustatnormal}), with
$H_n(X_j,X_k)=(n-1)^{-1}h^{d/2} (c_{jk}+c_{kj})$ where $c_{jk}
=a_{jk}-b_{jk}$, has been defined at the beginning of Section~\ref{s62}.
The asymptotic variance $v_1$ is the limit of the quantity $\frac{n^2}
2 \E[H_n(X_1,X_2)^2] $ asymptotically equivalent to
\begin{eqnarray*}
h^{d} \bigl(\E\bigl[c_{12}^2\bigr]+
\E[c_{12}c_{21}] \bigr).
\end{eqnarray*}
To compute this easily, we introduce the function $\xi
_i(x)=K_h(x-X_i)-f_h(x)$. First, use some algebra to obtain the formula
$b_{12} = \psi_{2}(X_1)\xi_2(X_1)-\int\psi_{2}(x)\xi_2(x)f(x)\,dx$, then
it follows that
\begin{eqnarray*}
c_{12}&=& a_{12}-b_{12}
\\
&=&\int\psi_3(x) \xi_1(x) \xi_2(x) f(x)\,dx -
\psi_{2}(X_1)\xi_2(X_1)+\int
\psi_{2}(x)\xi_2(x)f(x)\,dx
\\
&=&\int\bigl(\psi_3(x)f(x)\xi_2(x)-
\psi_2(X_1)\xi_2(X_1)
\bigr)K_h(x-X_1)\,dx
\\
&&{}+\int\psi_3(x)f(x)\xi_2(x) \bigl(f_h(x)-f(x)
\bigr)\,dx
\\
&=&\int\bigl(\psi_2(x)\xi_2(x)- \psi_2(X_1)
\xi_2(X_1) \bigr)K_h(x-X_1)\,dx
\\
&&{}+\int\psi_3(x)\xi_2(x) \bigl(f(x)-f_h(x)
\bigr)K_h(x-X_1)\,dx +\int\psi_3(x)f(x)
\xi_2(x) \bigl(f_h(x)-f(x) \bigr)\,dx
\\
& =& \int\bigl(\psi_2(x)K_h(x-X_2)-
\psi_2(X_1)K_h(X_1-X_2)
\bigr)K_h(x-X_1)\,dx
\\
&&{} +\int\bigl(- \bigl(\psi_2(x)f_h(x)-
\psi_2(X_1)f_h(X_1) \bigr)+
\psi_3(x)\xi_2(x) \bigl(f(x)-f_h(x) \bigr)
\bigr) K_h(x-X_1)\,dx
\\
&&{}+\int\psi_3(x)f(x)\xi_2(x) \bigl(f_h(x)-f(x)
\bigr)\,dx.
\end{eqnarray*}
Because $K_h$ integrates to $1$, it is not hard to see that the last
two terms in the previous equation will be negligible in the
computation of $v_1$. As a consequence, $ h^{d}\E[c_{12}^2]$ has the
same limit as
\begin{eqnarray*}
&& h^{d}\iint\biggl(\int\bigl(\psi_2(x)K_h(x-z)-
\psi_2(y)K_h(y-z) \bigr) K_h(x-y)\,dx
\biggr)^2f(y)f(z)\,dy\,dz
\\
&&\quad = \iint\biggl(\int\bigl(\psi_2(y+hu)K(u+v)-
\psi_2(y)K(v) \bigr) K(u)\,du \biggr)^2f(y)f(y-hv)\,dy\,dv
\\
&&\quad = V_K \int\psi_2(y)^2
f(y)^2\,dy +o(1)
\end{eqnarray*}
with $V_K = \int(\int(K(u+v)-K(v))K(u)\,du )^2\,dv$ and where the first
equality follows from a change of variables and the last representation
follows from the Lebesgue dominated theorem. Following the same steps
as previously, we obtain an similar expression for $h^d\E
[c_{12}c_{21}]$ and then we get
\begin{eqnarray*}
v_1 = 2V_K \int\varphi(y)^2f(y)^{-2}
\,dy.
\end{eqnarray*}
It remains to check the conditions of Theorem~\ref{ustatnormal}.
Clearly, the computation of $v_1$ provides that $\E
[H_n(X_1,X_2)^2]\approx n^{-2}$.
We obtain similarly that
$\E[H_n(X_1,X_2)^4]= O( n^{-5} h^{-d})$
and $ \E[G_n(X_1,X_2)^2]= O( n^{-5} h^{d})$ which implies the
conditions of the theorem.
\end{pf*}

\subsubsection{Proof of Theorem \texorpdfstring{\protect\ref{tclphinonregular}}{3}}

By (B1) and Lemma~\ref{holderremark}, there exists $M_2>0$
such that $\varphi$ is $\EuScript H (\min(s,1/2),M_2)$. Then we can
apply Lemma~\ref{decompustatmart} and by assumption on $h$, we obtain that
\begin{eqnarray*}
\bigl(nh^{-1}\bigr)^{1/2}(B_n+U_n+R_4+R_5)
&=& O_\P\bigl( n^{1/2}h^{r-1/2}+n^{-1/2}
h^{-(d+1)/2} +n^{-1}h^{-(3d+1)/2}\bigr)
\\
& =& o_{\P}(1).
\end{eqnarray*}
Since $M_n $ is a sum of independent variables with zero-mean, we can
apply the central limit theorem by checking the Lindeberg condition
(see, e.g., \cite{hall1980}, Chapter~3). Now we only have to
compute the asymptotic variance $v_2$ defined as the limit of
\begin{eqnarray*}
\var\bigl(h^{-1/2} \bigl(\psi_1(X_1) -
\widetilde\psi(X_1)\bigr) \bigr) = h^{-1}\E\bigl[\bigl(
\psi_1(X_1) - \widetilde\psi(X_1)
\bigr)^2\bigr]-h^{-1}\E\bigl[\psi_1(X_1)
- \widetilde\psi(X_1)\bigr]^2.
\end{eqnarray*}
On the one hand, by equations~(\ref{psitildeh}) and~(\ref{norml1psih}),
we have for some constant $C$
\begin{eqnarray*}
\bigl\llVert h^{-1/2}\bigl(\psi_{1h}(X_1) -
\widetilde\psi(X_1)\bigr)\bigr\rrVert_2&\leq&  C
h^{r-1/2},
\\
h^{-1/2} \bigl\llvert\E\bigl[\psi_1(X_1) -
\psi_{1h}(X_1)\bigr]\bigr\rrvert&\leq&  C
h^{r-1/2},
\end{eqnarray*}
as a consequence, we get
\begin{eqnarray*}
\var\bigl(h^{-1/2} \bigl(\psi_1(X_1) -
\widetilde\psi(X_1)\bigr) \bigr)= h^{-1}\bigl\llVert
\psi_1(X_1) - \psi_{1h}(X_1)
\bigr\rrVert_2^2 +o(1).
\end{eqnarray*}
On the other hand, for every $x\in Q$, we have
\begin{eqnarray*}
\psi_{1}(x)-\psi_{1h}(x) &=&\int_{{Q}^c}
\bigl(\psi_{1}(x)-\psi_{1}(y) \bigr)K_h(x-y)\,dy+
\int_{Q} \bigl(\psi_{1}(x)-\psi_{1}(y)
\bigr)K_h(x-y)\,dy
\\
&=& \psi_{1}(x)\int_{Q^c}K_h(x-y)\,dy
+ \int_{Q} \bigl(\psi_{1}(x)-
\psi_{1}(y) \bigr)K_h(x-y)\,dy,
\end{eqnarray*}
where $Q^c$ stands for the complement of the set $Q$ in $\R^d$. Because
$\psi_{1}$ is Nikolski with regularity $\min(s,r)$ inside $Q$, we use
equation~(\ref{norml2psih}) of Lemma~\ref{bochner} to show\vspace*{1pt} that the
$L_2$-norm of the right-hand side term is of order $h^{\min(s,r)}$. Clearly,
since $\min(s, r)>1/2$ we have
\begin{eqnarray*}
\var\bigl(h^{-1/2} \bigl(\psi_1(X_1) -
\widetilde\psi(X_1)\bigr) \bigr)= h^{-1}\biggl\llVert
\psi_{1}(X_1)\int_{Q^c}K_h(X_1-y)\,dy
\biggr\rrVert_2 +o(1)
\end{eqnarray*}
and it remains to apply Lemma~\ref{lemmageometry} to derive the stated limit.

\subsection{Proof of the Theorem \texorpdfstring{\protect\ref{thelemma2}}{4}}
By equation~(\ref{modelad}), we are interested in the asymptotic law of
\begin{eqnarray*}
n^{-1/2} \sum_{i=1}^n
\frac{\sigma(X_i) \psi(X_i)}{\w f_i} e_i + n^{-1/2} \Biggl(\sum
_{i=1}^n \frac{g(X_i) \psi(X_i)}{\w f_i}-\int g(x) \psi(x)\,d x
\Biggr).
\end{eqnarray*}
By Lemma~\ref{thelemma}, the right-hand side term goes to $0$ in
probability. For the other term, we use the decomposition $A_1 +A _2$, with
\begin{eqnarray*}
A_1 = n^{-1/2} \sum_{i=1}^n
\frac{\sigma\psi(X_i)}{ f(X_i)} e_i\quad\mbox{and}\quad A_2
=n^{-1/2} \sum_{i=1}^n
\frac{\sigma\psi(X_i)( f(X_i)-\w f(X_i))}{\w f_i f(X_i)} e_i,
\end{eqnarray*}
where $\sigma\psi(X_i) =\sigma(X_i) \psi(X_i)$. We define $\mathcal F
$ as the $\sigma$-field generated by the set of random variables $\{
X_1,X_2,\ldots\}$. We get
\begin{eqnarray*}
\E\bigl[A_2^2 \mid\mathcal F\bigr]= n^{-1}
\sum_{i=1}^n \frac{\sigma\psi(X_i)^2( f(X_i)-\w f_i)^2}{\w f_i^2 f(X_i)^2},
\end{eqnarray*}
then, one has
\begin{eqnarray*}
\E\bigl[A_2^2 \mid\mathcal F\bigr]\leq
\Bigl(b^2\inf_{i} \w f_i^2
\Bigr)^{-1} \llVert\sigma\psi\rrVert_{\infty}^2
n^{-1}\sum_{i=1}^n \bigl(
f(X_i)-\w f_i\bigr)^2.
\end{eqnarray*}
For the term on the left, since $\sigma\psi$ has support $Q$ we can
use~(\ref{unifb}), that is for $n$ large enough, it is bounded. For the
right-hand side term, it follows that
\begin{eqnarray*}
n^{-1}\sum_{i=1}^n \bigl(
f(X_i)-\w f_i\bigr)^2\leq
2n^{-1}\sum_{i=1}^n \bigl(
f(X_i)-f_h(X_i)\bigr)^2+2n^{-1}
\sum_{i=1}^n \bigl(f_h(X_i)-
\w f_i\bigr)^2,
\end{eqnarray*}
and then using equation~(\ref{bochner2}) in Lemma~\ref{bochner} and
(\ref{rose}) for $p=2$ we provide the bound
%
\begin{eqnarray}
\label{sansnom} \Biggl\llVert n^{-1}\sum_{i=1}^n
\bigl( f(X_i)-\w f_i\bigr)^2\Biggr\rrVert
_1\leq  C\bigl(h^{2r}+n^{-1}h^{-d}
\bigr)
\end{eqnarray}
for some $C>0$.
Therefore, we have shown that
$\E[A_2^2 \mid\mathcal F]\r0 $ in probability. Since for any
$\varepsilon>0$, $\P( \llvert A_2\rrvert >\varepsilon\mid\mathcal F)
\leq \varepsilon^{-2} \E[A_2^2 \mid\mathcal F]$, it remains to note
that the sequence $\P(\llvert A_2\rrvert >\varepsilon\mid\mathcal F)$
is uniformly integrable to apply the Lebesgue domination theorem to get
\begin{eqnarray*}
\P(A_2>\varepsilon)\lr0.
\end{eqnarray*}
To conclude, we apply the central limit theorem to $A_1$ and the
statement follows.

\subsection{Some lemmas}

\subsubsection{Inequalities}

\begin{lemma}\label{bochner} For any function $g:\R^d \r\R$, recall
that $ g_h(x) = \int g(x-hu)K(u)\,du$. Under assumptions \textup{(A1)},
\textup{(A2)} and \textup{(A4)}, it holds that
%
\begin{eqnarray}
\label{bochner2} \llVert f_h - f\rrVert_\infty&\leq&
C_Kh^r\bigl\llVert f^{(r)}\bigr\rrVert
_\infty,
\\
\bigl\llvert\E\bigl[\varphi(X_1) -\varphi_{h}(X_1)
\bigr]\bigr\rrvert&\leq&  C_K h^r\bigl\llVert
f^{(r)}\bigr\rrVert_\infty\int\bigl\llvert\varphi(x)\bigr
\rrvert \,dx, \label{norml1psih}
\\
\bigl\llVert\varphi_h(X_1) -
\varphi(X_1)\bigr\rrVert_2 &\leq&  C_K M
h^{s}, \label{norml2psih}
\end{eqnarray}
where $C_K$ is a positive constant that depends $K$ only.
\end{lemma}

\begin{pf}
We start by proving~(\ref{norml1psih}) and~(\ref{norml2psih}) assuming
that~(\ref{bochner2}) holds. For the mean: using Fubini's theorem, we have
\begin{eqnarray*}
\E\bigl[\varphi(X_1) -\varphi_{h}(X_1)
\bigr] &=&\int\bigl(\varphi(x)-\varphi_h(x) \bigr)f(x) \,dx
\\
&=&\int\varphi(x)f(x)-\varphi(x) f_{h}(x) \,dx,
\end{eqnarray*}
hence
\begin{eqnarray*}
\bigl\llvert\E\bigl[\varphi(X_1) -\varphi_{h}(X_1)
\bigr]\bigr\rrvert&\leq& \bigl\llVert f(x)-f_{h}(x)\bigr\rrVert
_\infty\int\bigl\llvert\varphi(x)\bigr\rrvert \,dx,
\end{eqnarray*}
which by~(\ref{bochner2}) gives
\begin{eqnarray*}
\bigl\llvert\E\bigl[\varphi(X_1) -\varphi_{h}(X_1)
\bigr]\bigr\rrvert&\leq&  C_Kh^r\bigl\llVert
f^{(r)}\bigr\rrVert_\infty\int\bigl\llvert\varphi(x)\bigr
\rrvert \,dx.
\end{eqnarray*}
This is~(\ref{norml1psih}). We turn now to~(\ref{norml2psih}):
%
\begin{eqnarray}
\label{dfg} \E\bigl[\bigl(\varphi_h(X_1) -
\varphi(X_1)\bigr)^2\bigr] &=& \int\biggl(\int\bigl(
\varphi(x-hu)-\varphi(x) \bigr)K(u)\,du \biggr)^2f(x) \,dx.
\end{eqnarray}
We now use the Taylor formula with Lagrange remainder applied to
$g(t)=\varphi(x-tu)$ with
order $k$ equal to the largest integer smaller than $s$:
\begin{eqnarray*}
\varphi(x-hu)-\varphi(x)&=&\sum_{j=1}^{k-1}
\frac{h^j}{j!}g^{(j)}(0)+ \int_0^{h}g^{(k)}(t)
\frac{(h-t)^{k-1}}{(n-1)!}\,dt
\\
&=&\sum_{j=1}^k\frac{h^j}{j!}g^{(j)}(0)+
\int_0^{h}\bigl(g^{(k)}(t)-g^{(k)}(0)
\bigr)\frac{(h-t)^{k-1}}{(n-1)!}\,dt.
\end{eqnarray*}
The first term is a polynomial in $u$ which will vanish
after insertion in~(\ref{dfg}) because $K$ is orthogonal the first
non-constant polynomial of degree
$\leq  r$. The second term is bounded as
\begin{eqnarray*}
\biggl\llvert\int_0^{h}\bigl(g^{(k)}(t)-g^{(k)}(0)
\bigr)\frac{(h-t)^{k-1}}{(k-1)!}\,dt\biggr\rrvert\leq \llvert
u\rrvert
^kh^{k-1}\int_0^{h}\bigl
\llvert\varphi^{(k)}(x-tu)-\varphi^{(k)}(x)\bigr\rrvert \,dt.
\end{eqnarray*}
Hence,
%
\begin{eqnarray}\label{taylor}
&& \biggl\llvert\int\bigl(\varphi(x-hu)-\varphi(x) \bigr)K(u)\,du
\biggr\rrvert
\nonumber\\[-8pt]\\[-8pt]\nonumber
&&\quad \leq  h^{k-1}\int_0^{h}
\int\bigl\llvert\varphi^{(k)}(x-tu)-\varphi^{(k)}(x)\bigr
\rrvert\llvert u\rrvert^k K(u)\,du \,dt
\end{eqnarray}
and by the generalized Minkowski inequality (\cite{folland1999} page~194)\footnote
{For any non-negative measurable function $g(\cdot,\cdot)$ on $\mathbb R^{k+d}$,
\begin{eqnarray*}
\biggl(\int\biggl(\int g(y,x)\,dy \biggr)^2 \,dx \biggr)^{1/2}&
\leq& \int\biggl(\int g(y,x)^2\,dx \biggr)^{1/2} \,dy.
\end{eqnarray*}
}
\begin{eqnarray*}
\bigl\llVert(\varphi_h -\varphi) (X_1)\bigr\rrVert
_2 &\leq&  h^{k-1}\int\biggl(\int\bigl\llvert
\varphi^{(k)}(x-tu)-\varphi^{(k)}(x)\bigr\rrvert
^2u^{2k}K(u)^2 1_{0\leq  t\leq  h}f(x)\,dx
\biggr)^{1/2} \,du\,dt
\nonumber
\\
&\leq&  M h^{k-1}\int\bigl(\llvert tu\rrvert^{2\alpha}
\llvert u\rrvert^{2k}K(u)^2 \bigr)^{1/2}
1_{0\leq  t\leq  h} \,du\,dt
\nonumber
\\
&=& M(1+\alpha)^{-1} h^{k+\alpha}\int\bigl(\llvert u\rrvert
^{2\alpha+2k}K(u)^2 \bigr)^{1/2} \,du
\nonumber
.
\end{eqnarray*}
This implies~(\ref{norml2psih}).
Concerning~(\ref{bochner2}), we use~(\ref{taylor}) with $f$ and $k=r$
to get that
\begin{eqnarray*}
\bigl\llvert f_h(x)-f(x)\bigr\rrvert&\leq&  h^{r-1}
\int_0^{h}\int\bigl\llvert f^{(r)}(x+tu)
\bigr\rrvert\llvert u\rrvert^r K(u)\,du \,dt,
\end{eqnarray*}
the latter is bounded by a constant times $h^r$.
\end{pf}

The following lemma gives some bounds on the conditional moments of $\xi
_{12}$ that are useful in the proof of Theorem~\ref{thelemma}.

\begin{lemma}\label{moments}
Let $\xi_{ij}=K_{ij} - f_h(X_i)$, under \textup{(A1)} and \textup{(A2)}
%
\begin{eqnarray}
\E[\xi_{12}\mid X_1]&=&0,\label{eq1}
\\
\E\bigl[\llvert\xi_{12}\rrvert^p\mid X_1
\bigr]&\leq& 2^p\E\bigl[\llvert K_{12}\rrvert
^p\mid X_1 \bigr]\leq  C_1
h^{-(p-1)d},\label{eq2}
\\
\bigl\llvert\E[\xi_{13}\xi_{23}\mid X_1,X_2]
\bigr\rrvert&\leq& \llVert f\rrVert_\infty\bigl(h^{-d}
\widetilde K\bigl(h^{-1}(X_2-X_1)\bigr)+\llVert
f\rrVert_{\infty} \bigr)\label{eq3},
\end{eqnarray}
with $\widetilde K(x)=\int\llvert K(x-y)K(y)\rrvert \,dy$ and $C_1>0$.\vadjust{\goodbreak}
\end{lemma}

\begin{pf}
The first equation is trivial. For the second equation, the triangular
inequality and the Jensen inequality provide
\begin{eqnarray*}
\E\bigl[\llvert\xi_{12}\rrvert^p\mid X_1
\bigr]\leq 2^p \E\bigl[\llvert K_{12}\rrvert
^p\mid X_1\bigr]=2^p h^{-(p-1)d}\int
\bigl\llvert K(u)\bigr\rrvert^pf(X_1-h u) \,d x,
\end{eqnarray*}
and the third one is derived by
\begin{eqnarray*}
\bigl\llvert\E[\xi_{13}\xi_{23}\mid X_1,X_2]
\bigr\rrvert&=&\bigl\llvert\E[\xi_{13}K_{23}\mid
X_1,X_2]\bigr\rrvert
\\
&=&\biggl\llvert\int\bigl(K_h(X_1-x)-f_h(X_1)
\bigr)K_h(X_2-x)f(x)\,dx\biggr\rrvert
\\
&=&\biggl\llvert\int\bigl(K_h(X_1-X_2+hu)-f_h(X_1)
\bigr)K(u)f(X_2-hu)\,du\biggr\rrvert
\\
&\leq& \llVert f\rrVert_\infty\bigl(h^{-d}\widetilde K
\bigl(h^{-1}(X_2-X_1)\bigr)+\llVert f\rrVert
_{\infty}\bigr).
\end{eqnarray*}\upqed
\end{pf}

The Efron--Stein inequality helps to bound the $L_2$ moments of
estimators. For the proof, we refer to the original paper \cite
{efron1981} but also to \cite{boucheron2004}.

\begin{theorem}[(Efron--Stein inequality)]\label{efronstein}
Let $X_1,\dots, X_n$ be an i.i.d. sequence, $X_1'$ be an independent
copy of $X_1$
and $f$ be a symmetric function of $n$ variables, then
\begin{eqnarray*}
\var\bigl(f(X_1,\dots,X_n)\bigr)\leq
\frac{n} 2\E\bigl[\bigl(f(X_1,\dots, X_n)-f
\bigl(X'_1,X_2,\dots, X_n\bigr)
\bigr)^2 \bigr].
\end{eqnarray*}
\end{theorem}

\subsubsection{Measure results}

\begin{lemma}\label{holderremark}
Let\vspace*{1pt} $s>0$ and $M_1>0$, suppose that the support of $\varphi$ is a
convex body $Q$ and that $\varphi$
is $\EuScript H(s,M_1)$ on $Q$, then there exists $M >0$ such that
$\varphi$ is $\EuScript H(\min(s,1/2),M)$ on $\R^d$.
\end{lemma}
\begin{pf} We have
\begin{eqnarray*}
&& \int \bigl\llvert\varphi(x+u)-\varphi(x)\bigr\rrvert^2 \,dx
\\
&&\quad = \int_{\{ x\in Q,x+u\in Q\}} \bigl\llvert\varphi(x+u)-\varphi(x)\bigr
\rrvert^2 \,dx + \int_{\{ x\notin Q,x+u\in Q\}} \varphi(x+u)^2
\,dx
\\
&&\qquad{} + \int_{\{ x\in Q,x+u\notin Q\}} \varphi(x)^2 \,dx
\\
&&\quad = \int_{\{ x\in Q,x+u\in Q\}} \bigl\llvert\varphi(x+u)-\varphi(x)\bigr
\rrvert^2 \,dx + \int_Q \varphi(x)^2
( 1 _{\{x-u\notin Q\}}+ 1 _{\{x+u\notin Q\}}) \,dx
\\
&&\quad \leq \int_{\{ x\in Q,x+u\in Q\}} \bigl\llvert\varphi(x+u)-\varphi(x)
\bigr\rrvert^2 \,dx + \llVert\varphi\rrVert_{\infty}^2
\int1_{\{\operatorname{dist}(x,\partial Q)\leq \llvert u\rrvert
\}} \,dx
\\
&&\quad \leq  M_1\llvert u\rrvert^{2s }+ \llVert\varphi
\rrVert_{\infty}^2 \xi_{d-1}(Q) \llvert u\rrvert,
\end{eqnarray*}
where $\xi_{d-1}(S)$ is called a Quermassintegrale of Minkowski and
$\operatorname{dist}$ stands for the
Euclidean distance in $\R^d $.
The last inequality follows from the fact that $\varphi$ is $\EuScript H(s,M_1)$
on $Q$ and by the Steiner's formula stated, for instance, in \cite
{federer1969}, Theorem 3.2.35, page 271.
\end{pf}

\begin{lemma}\label{lemmageometry}
Under the assumption \textup{(A4)}, if $Q$ is a compact set with
$\mathcal C^2$ boundary and
$\psi$ is continuous
\begin{eqnarray*}
&&\lim_{h\rightarrow0}h^{-1}\int_Q
\biggl(\int_{Q^c} K_h(x-y)\,dy
\biggr)^2\psi(x)\,dx =\int_{\partial Q} L_Q(x)
\psi(x) \,d\mathcal H^{d-1}(x),
\end{eqnarray*}
where
\begin{eqnarray*}
L_Q(x)=\iint\min\bigl(\bigl\langle z,u(x)\bigr\rangle,\bigl\langle
z',u(x)\bigr\rangle\bigr)_+K(z)K\bigl(z'
\bigr)\,dz\,dz',
\end{eqnarray*}
and $\mathcal H^{d-1}$ stands for the $(d-1)$-dimensional Hausdorff
measure, $u(x)$ is the normal outer vector of $Q$ at the point $x$.
\end{lemma}

\begin{pf} Let us start with an estimate of the integral over $Q^c$
having a simpler dependency w.r.t.~$h$.
We define the function
\begin{eqnarray*}
\tau(x)=(1_{x\in Q}-1_{x\notin Q})\operatorname{dist}(x,\partial Q).
\end{eqnarray*}
This function is $\mathcal C^2$ in the neighborhood of $\partial Q$ and
its gradient $-u(x)$
is, for $x\in\partial Q$, the normal inner vector (since $\partial Q$
is $\mathcal C^2$,
using a local parametrization of $Q$, we are reduced to the case where
$\partial Q$
is a piece of hyperplane).
Then
\begin{eqnarray*}
\int_{Q^c} K_h(x-y)\,dy &=&\int
1_{x+hz\in Q^c}K(z)\,dz
\\
&=&\int1_{\tau(x+hz)\leq 0}K(z)\,dz
\\
&=&\int1_{\tau(x)- h\langle z,u(x)\rangle\leq  ah^2}K(z)\,dz,
\end{eqnarray*}
where $a$ actually depends on $x$ and $z$ but is smaller than a
constant related to
the curvature of~$\partial Q$. Hence,
\begin{eqnarray*}
&& \biggl\llvert\int_{Q^c} K_h(x-y)\,dy-\int
1_{\tau(x)- h\langle z,u(x)\rangle\leq 0}K(z)\,dz\biggr\rrvert
\\
&&\quad \le\int\llvert1_{\tau(x)- h\langle z,u(x)\rangle\leq
ah^2}-1_{\tau(x)- h\langle z,u(x)\rangle\le0}\rrvert K(z)\,dz
\\
&&\quad \le\int1_{\llvert \tau(x)- h\langle z,u(x)\rangle\rrvert \leq
\llvert a\rrvert h^2}K(z)\,dz
\\
&&\quad \leq  a_0h1_{\tau(x)\leq  m_0h}
\end{eqnarray*}
for some $a_0$ and $m_0$, because the integration domain is a band of width
$\llvert a\rrvert h$. Hence,
\begin{eqnarray*}
\int_{Q^c} K_h(x-y)\,dy & =& \int
1_{\tau(x)\leq  h\langle z,u(x)\rangle}K(z)\,dz+a_1(x)h1_{\tau
(x)\leq  m_0h},
\end{eqnarray*}
where $a_1$ is bounded. Since the second
term has a $O(h^2)$ integral over $Q$, its contribution in the limit is
negligible,
and it suffices to prove that
\begin{eqnarray*}
\lim_{h\rightarrow0}h^{-1}\int_Q
\biggl(\int1_{\tau(x)\leq  h\langle z,u(x)\rangle}K(z)\,dz \biggr
)^2\psi(x)\,dx =\int
_{\partial Q}L_Q(x) \psi(x)\,d\mathcal H^{d-1}(x).
\end{eqnarray*}
By setting
\begin{eqnarray*}
\varphi(x,t)= \biggl(\int1_{0\leq  t\leq \langle z,u(x)\rangle
}K(z)\,dz \biggr)^2
\psi(x)1_{x\in Q},
\end{eqnarray*}
the latter equality can be rewritten as
\begin{eqnarray*}
\lim_{h\rightarrow0}h^{-1}\int\varphi\bigl(x,h^{-1}
\tau(x)\bigr)\,dx =\int_{\partial Q}L_Q(x)\psi(x) \,d\mathcal
H^{d-1}(x).
\end{eqnarray*}
From Proposition 3, page 118 of \cite{evans1992}, we have for any
integrable function $q$
and $f$ Lipschitz with $\operatorname{essinf}\llvert \nabla f\rrvert >0$:
\begin{eqnarray*}
\int_{f\geq 0} q(x)\,dx &=&\int_0^\infty
\biggl(\int_{f=s}\frac{q(x)}{\llvert \nabla f(x)\rrvert }\,d\mathcal
H^{d-1}(x) \biggr)\,ds
\end{eqnarray*}
hence, with $f(x)=h^{-1}\tau(x)$ and $q(x)=\varphi(x,h^{-1}\tau(x))$,
we obtain
\begin{eqnarray*}
h^{-1}\int\varphi\bigl(x,h^{-1}\tau(x)\bigr)\,dx &=&\int
_0^\infty\biggl(\int_{\tau=hs}
\varphi(x,s) \,d\mathcal H^{d-1}(x) \biggr)\,ds.
\end{eqnarray*}
Letting $h\rightarrow0$, we get
\begin{eqnarray*}
\lim_{h\rightarrow0}h^{-1}\int\varphi\bigl(x,h^{-1}
\tau(x)\bigr)\,dx &=&\int_0^\infty\biggl(\int
_{\partial Q}\varphi(x,s)\,d\mathcal H^{d-1}(x) \biggr)\,ds
\\
&=&\int_{\partial Q} \biggl(\int_0^\infty
\varphi(x,s)\,ds \biggr)\,d\mathcal H^{d-1}(x).
\end{eqnarray*}
We can write $\int\varphi(x,s)\,ds$ as
\begin{eqnarray*}
\int_{0}^\infty\varphi(x,s)\,ds &=&\psi(x)\iiint
1_{0\leq  s\leq \langle z,u(x)\rangle}1_{0\leq  s\leq
\langle z',u(x)\rangle}K(z)K\bigl(z'
\bigr)\,dz\,dz'\,ds
\\
&=&\psi(x)\iint\min\bigl(\bigl\langle z,u(x)\bigr\rangle,\bigl\langle
z',u(x)\bigr\rangle\bigr)_+K(z)K\bigl(z'
\bigr)\,dz\,dz'.
\end{eqnarray*}\upqed
\end{pf}

\subsubsection{Weak convergence for degenerate $U$-statistics}

\begin{theorem}[(Hall (1984), \cite{hall1984})]\label{ustatnormal}
Let $H_n:\R^d\times\R^d \r\R$, with $H_n$ symmetric, assume that $\E
[H_n(X_1,X_2)\mid X_1]=0$ and $\E[H_n(X_1,X_2)^2]<+\infty$. If
\begin{eqnarray*}
\frac{\E[G_n(X_1,X_2)^2]+n^{-1} \E[H_n(X_1,X_2)^4]}{\E
[H_n(X_1,X_2)^2]^2} \stackrel{n\r
+\infty} {\lr}
0,
\end{eqnarray*}
with $G_n(x,y)=\E[H_n(X_1,x)H_n(X_1,y)]$, then $\sum_{j<k} H(X_j,X_k) $
is asymptotically normally distributed with zero mean and variance
given by $\frac{n^2} 2 \E[H(X_1,X_2)^2 ] $.
\end{theorem}


\section*{Acknowledgements}
The authors would like to thank C\' eline Vial for helpful comments and
advice on a latter version of this article.

Research supported by the Fonds de la Recherche Scientifique (FNRS) A4/5 FC 2779/2014--2017 No. 22342320.


%

\printhistory
\end{document}